\documentclass[10pt, reqno]{amsart}

\usepackage{amssymb}
\usepackage{hyperref} 
\usepackage{mathrsfs,bbold}
\usepackage{amsmath}
\usepackage{amsthm}
\usepackage{enumerate}
\usepackage{dsfont}
\usepackage{color}
\usepackage[a4paper]{geometry}
\usepackage[utf8]{inputenc}
\usepackage{stmaryrd}
\usepackage{titlesec}
\usepackage{mathabx}
\usepackage{appendix}
\usepackage{todonotes, fancyhdr, wasysym }

\usepackage{tikz-cd}
\usetikzlibrary{calc}

\usepackage[all]{xy}
\let\oldtocsection=\tocsection
\let\oldtocsubsection=\tocsubsection
\let\oldtocsubsubsection=\tocsubsubsection

\renewcommand{\tocsection}[2]{\hspace{0em}\oldtocsection{#1}{#2}}
\renewcommand{\tocsubsection}[2]{\hspace{1em}\oldtocsubsection{#1}{#2}}
\renewcommand{\tocsubsubsection}[2]{\hspace{2em}\oldtocsubsubsection{#1}{#2}}

\numberwithin{equation}{section}

\newcommand{\ssk}{\smallskip}

\renewcommand{\epsilon}{\varepsilon}

\newcommand\bbE{\mathbb{E}}

\newcommand\bbN{\mathbb{N}}
\newcommand\bbP{\mathbb{P}}
\newcommand\bbQ{\mathbb{Q}}
\newcommand\bbR{\mathbb{R}}

\newcommand{\mcC}{\mathcal{C}} 
\newcommand{\mcD}{\mathcal{D}}
\newcommand{\mcE}{\mathcal{E}}
\newcommand{\mcF}{\mathcal{F}}

\newcommand{\mcL}{\mathcal{L}}

\newcommand{\mcP}{\mathcal{P}}

\newcommand\mcT{\mathcal T}
\newcommand\mcV{\mathcal V}

\newenvironment{Dem}[1][\unskip]{%
    \begin{list}{\hspace{1.15cm}{\sf \textbf{Proof #1 --}}}{%
        \setlength{\topsep}{0pt}%
        \setlength{\leftmargin}{0pt}%
        \setlength{\rightmargin}{0pt}%
        \setlength{\listparindent}{0pt}%
        \setlength{\itemindent}{0pt}%
        \setlength{\parsep}{0pt}%
        \addtolength{\leftmargin}{0pt}
        \addtolength{\rightmargin}{0pt}%
    } \item }{\hfill $\rhd$\end{list}\smallskip}

\titleformat{\section}[block]
{\filcenter\normalfont\sffamily\bfseries\Large}
{{\hspace{-0.7cm}}\thesection \hspace{0.2em} --\vspace{0.15cm}}{0.5em}{}

\titleformat{\subsection}[runin]
{\filcenter\normalfont\sffamily\bfseries\large}  						  
{\hspace{0cm}\thesubsection \hspace{0.5em}--\vspace{0cm}}{.5em}{}  
\titlespacing{\subsection}{-0pc}{1.5ex plus .1ex minus .2ex}{0pc}

\titleformat{\subsubsection}[runin]
{\filcenter\normalfont\bfseries\sffamily}
{{\hspace{-0cm}}{\thesubsubsection} \hspace{-0.1em} \hspace{-0.2cm}}{0.5em}{}  

\numberwithin{subsection}{section}
\numberwithin{subsubsection}{subsection}

\newtheoremstyle{mystyle}
{3pt}               
{3pt}               
{\it }                      
{}                      
{\sffamily\bfseries}             
{}                      
{0.5em}                 
{{#2.} #1 { --}}   

\theoremstyle{mystyle}

\newtheorem{thm}{Theorem}[section]
	\numberwithin{thm}{section}
	
\newtheorem{cor}[thm]{\hspace{-0.15cm}  {Corollary} }
\newtheorem{lem}[thm]{\hspace{-0.2cm}  {Lemma} }
\newtheorem*{defn*} {Definition}
\newtheorem*{prop*} {Proposition}
\newtheorem*{lem*} {Lemma}
\newtheorem*{cor*} {Corollary}

\newtheoremstyle{mystyle2}
{3pt}               
{3pt}               
{\it }                      
{}                      
{\sffamily\bfseries}             
{}                      
{0.5em}                 
{\llap{#2 }#1{\hspace{0.2cm}--}}
\theoremstyle{mystyle2}

\newtheorem*{definition*}{Definition}
\newtheorem*{lemma*}{lem}
\newtheorem*{thm*}{thm}
\numberwithin{equation}{section} 
\newtheorem*{Remark*}{Remark}

\newcommand\restr[2]{{
  \left.\kern-\nulldelimiterspace 
  #1 
  \vphantom{\big|} 
  \right|_{#2} 
  }}

\makeatletter
\newcommand*{\defeq}{\mathrel{\rlap{%
                     \raisebox{0.25ex}{$\m@th\cdot$}}%
                     \raisebox{-0.25ex}{$\m@th\cdot$}}%
                     =}
\makeatother

\makeatletter
\newcommand*{\eqdef}{=\mathrel{\rlap{%
                     \raisebox{0.25ex}{$\m@th\cdot$}}%
                     \raisebox{-0.25ex}{$\m@th\cdot$}}%
                     }
\makeatother

\usepackage{tikz-cd}
\usetikzlibrary{calc}
\usetikzlibrary{shapes.misc}
\usetikzlibrary{shapes.symbols}
\usetikzlibrary{snakes}
\usetikzlibrary{decorations}
\usetikzlibrary{decorations.markings}

\makeatletter
\pgfdeclareshape{crosscircle}
{
  \inheritsavedanchors[from=circle] 
  \inheritanchorborder[from=circle]
  \inheritanchor[from=circle]{north}
  \inheritanchor[from=circle]{north west}
  \inheritanchor[from=circle]{north east}
  \inheritanchor[from=circle]{center}
  \inheritanchor[from=circle]{west}
  \inheritanchor[from=circle]{east}
  \inheritanchor[from=circle]{mid}
  \inheritanchor[from=circle]{mid west}
  \inheritanchor[from=circle]{mid east}
  \inheritanchor[from=circle]{base}
  \inheritanchor[from=circle]{base west}
  \inheritanchor[from=circle]{base east}
  \inheritanchor[from=circle]{south}
  \inheritanchor[from=circle]{south west}
  \inheritanchor[from=circle]{south east}
  \inheritbackgroundpath[from=circle]
  \foregroundpath{
    \centerpoint%
    \pgf@xc=\pgf@x%
    \pgf@yc=\pgf@y%
    \pgfutil@tempdima=\radius%
    \pgfmathsetlength{\pgf@xb}{\pgfkeysvalueof{/pgf/outer xsep}}%
    \pgfmathsetlength{\pgf@yb}{\pgfkeysvalueof{/pgf/outer ysep}}%
    \ifdim\pgf@xb<\pgf@yb%
      \advance\pgfutil@tempdima by-\pgf@yb%
    \else%
      \advance\pgfutil@tempdima by-\pgf@xb%
    \fi%
    \pgfpathmoveto{\pgfpointadd{\pgfqpoint{\pgf@xc}{\pgf@yc}}{\pgfqpoint{-0.707107\pgfutil@tempdima}{-0.707107\pgfutil@tempdima}}}
    \pgfpathlineto{\pgfpointadd{\pgfqpoint{\pgf@xc}{\pgf@yc}}{\pgfqpoint{0.707107\pgfutil@tempdima}{0.707107\pgfutil@tempdima}}}
    \pgfpathmoveto{\pgfpointadd{\pgfqpoint{\pgf@xc}{\pgf@yc}}{\pgfqpoint{-0.707107\pgfutil@tempdima}{0.707107\pgfutil@tempdima}}}
    \pgfpathlineto{\pgfpointadd{\pgfqpoint{\pgf@xc}{\pgf@yc}}{\pgfqpoint{0.707107\pgfutil@tempdima}{-0.707107\pgfutil@tempdima}}}
  }
}
\makeatother

\colorlet{symbols}{black}    
\colorlet{testcolor}{green!60!black}
\colorlet{supcolor}{red!60!black}

\tikzset{
	root/.style={circle, fill=testcolor!70, draw=testcolor, inner sep=1pt, minimum size=0.5mm},
	dot/.style={circle, draw=black, fill=black, inner sep=0pt, minimum size=0mm},
	noise/.style={circle, draw=black, fill=white, inner sep=0pt, minimum size=1mm},
	h/.style={circle, draw=black, fill=black, inner sep=0pt, minimum size=0.3mm},	
	var/.style={circle, fill=white, draw=purple, fill=purple, inner sep=0pt, minimum size=0.6mm},
	bdot/.style={circle, draw=black, fill=white, inner sep=0pt, minimum size=1mm},
	bluedot/.style={circle,draw=blue, fill=blue, inner sep=0pt, minimum size=2mm},
	noiseblue/.style={circle, fill=blue!20, draw=blue, inner sep=0pt, minimum size=1mm},
	dtestfcn/.style={ultra thick, densely dashed, testcolor,shorten >=1pt,shorten <=1pt,<-},
	testfcnx/.style={semithick,testcolor,shorten >=1pt,shorten <=1pt,<-, postaction={decorate,decoration={markings,mark=at position 0.6 with {\drawx}}}},
	testfcn/.style={semithick, testcolor, shorten >=1pt,shorten <=1pt,-},	
	K/.style= {semithick, shorten >=0pt,shorten <=0pt,-},
	kdashed/.style= {semithick,densely dashed,shorten >=1pt,shorten <=1pt,<->},
	DK/.style={thick, densely dotted, shorten >=0pt,shorten <=0pt},   
	ksup/.style={thick, supcolor, shorten >=1pt,shorten <=1pt},
	dots/.style={semithick,dotted,shorten >=1pt,shorten <=1pt},
	Deps/.style={semithick,draw=black!25,fill=black!25,shorten >=1pt,shorten <=1pt,->},
	kbase/.style={semithick,dotted,shorten >=1pt,shorten <=1pt,->},
	multx/.style={shorten >=1pt,shorten <=1pt,
		postaction={decorate,decoration={markings,mark=at position 0.5 with {\drawx}}}},
	kernelx/.style={semithick,shorten >=1pt,shorten <=1pt,->,
		postaction={decorate,decoration={markings,mark=at position 0.4 with {\drawx}}}},
	kernel1/.style={->,semithick,shorten >=1pt,shorten <=1pt,postaction={decorate,decoration={markings,mark=at position 0.45 with {\draw[--] (0,-0.1) -- (0,0.1);}}}},
	kernel2/.style={->,semithick,shorten >=1pt,shorten <=1pt,postaction={decorate,decoration={markings,mark=at position 0.45 with {\draw[--] (0.05,-0.1) -- (0.05,0.1);\draw[--] (-0.05,-0.1) -- (-0.05,0.1);}}}},
	kernelBig/.style={semithick,shorten >=1pt,shorten <=1pt,decorate, decoration={zigzag,amplitude=1.5pt,segment length = 3pt,pre length=2pt,post length=2pt}},
	rho/.style={dotted,semithick,shorten >=1pt,shorten <=1pt},
	renorm/.style={shape=circle,fill=white,inner sep=1pt},
	labl/.style={shape=rectangle,fill=white,inner sep=1pt},
	xi/.style={circle,fill=symbols!10,draw=symbols,inner sep=0pt,minimum size=1.2mm},
	xix/.style={crosscircle,fill=symbols!10,draw=symbols,inner sep=0pt,minimum size=1.2mm},
	xib/.style={circle,fill=symbols!10,draw=symbols,inner sep=0pt,minimum size=1.6mm},
	xibx/.style={crosscircle,fill=symbols!10,draw=symbols,inner sep=0pt,minimum size=1.6mm},
	not/.style={circle,fill=symbols,draw=symbols,inner sep=0pt,minimum size=0.5mm},
	>=stealth,
	graydot/.style={circle,fill=gray,inner sep=0pt, minimum size=1mm},
	zero/.style={circle,inner sep=0pt, minimum size=1mm, draw},
	kernelprimeeps/.style={densely dashed, semithick,shorten >=1pt,shorten <=1pt},
	smalldot/.style={circle,fill=black,draw=black, solid,inner sep=0pt,minimum size=0.5mm},
	}

\begin{document}

\begin{center}
{\LARGE\sffamily{\textbf{Uniqueness of the $\Phi^4_3$ measures on closed Riemannian $3$-manifolds}   \vspace{0.5cm}}}
\end{center}

\begin{center}
{\sf I. BAILLEUL\footnote{Univ Brest, CNRS, LMBA - UMR 6205, F- 29238 Brest, France. {\it E-mail}: ismael.bailleul@univ-brest. Partial support from the ANR via the ANR-22-CE40-0017 grant.} }
\end{center}

\vspace{1cm}

\begin{center}
\begin{minipage}{0.8\textwidth}
\renewcommand\baselinestretch{0.7} \scriptsize \textbf{\textsf{\noindent Abstract.}} We constructed in a previous work the $\Phi^4_3$ measures on compact boundaryless $3$-dimensional Riemannian manifolds as some invariant probability measures of some Markovian dynamics; a crucial feature to characterize the $\Phi^4_3$ measure. We prove in the present work that these dynamics have a unique invariant probability measures. This is done by using an explicit coupling by change of measure that does not require any a priori information on the support of the law of the solution to the dynamics. In addition, the coupling can be used to see that the semigroup generated by the dynamics satisfies a Harnack-type inequality, which entails that the semigroup has the strong Feller property.
\end{minipage}
\end{center}

\section{Introduction}
\label{SectionIntro}

Let $M$ stand for an arbitrary compact boundaryless Riemannian manifold of dimension $3$. Following the flow of research that grew out of the recent development of the domain of singular stochastic partial differential equations (PDEs), we constructed in our previous work \cite{BDFT} the $\Phi^4_3$ measure on $M$ as an invariant probability measure of some Markovian dynamics on the space $C^{-1/2-\epsilon}(M)$ of $(-\frac{1}{2}-\epsilon)$-H\"older/Besov distributions on $M$, for an arbitrary $\epsilon>0$. The dynamics is given by Parisi \& Wu's paradigm of stochastic quantization and takes the form
\begin{equation} \label{EqPhi43SPDE}
\partial_t u = (\Delta-1) u - u^3 + \xi,
\end{equation}
where $\xi$ stands for a spacetime white noise. When set on a discrete $3$-dimensional torus this PDE rewrites as a coupled system of stochastic differential equations whose invariant measure is unique and has a density with respect to the massive discrete Gaussian free field measure proportional to $\exp\big(-\frac{1}{4}\sum_i \phi_i^4\big)$. Its continuous counterpart is `the' $\Phi^4_3$ measure; it has density $\exp\big(-\frac{1}4{\int_M\phi^4}\big)$ with respect to the massive Gaussian free field measure on $M$. However this reference measure has support in the spaces $C^{-1/2-\epsilon}(M)$, for all $\epsilon>0$, and essentially no better. The fourth power of $\phi$ is thus almost surely ill-defined and a renormalization procedure is needed to construct such a measure from its density. In particular this makes the unique characterization of the $\Phi^4_3$ measure a non-trivial question. The stochastic quantization approach to the construction of the $\Phi^4_3$ measure postulates that Equation \eqref{EqPhi43SPDE} is well-posed for all times and that it defines a Markovian dynamics which has a unique invariant probability measure, defined as the $\Phi^4_3$ measure. This approach to the construction of the $\Phi^4_3$ measure does not avoid the need of a renormalization process. Indeed spacetime white noise in space dimension $3$ has H\"older parabolic regularity $-5/2-\epsilon$, for all $\epsilon>0$, and no better, so a solution to Equation \eqref{EqPhi43SPDE} has at best parabolic H\"older regularity $-1/2-\epsilon$, and the quantity $u^3$ is ill-defined. This problem is what makes Equation \eqref{EqPhi43SPDE} a {\it singular} stochastic PDE. Its proper formulation requires a priori the use of an ad hoc setting such as regularity structures \cite{Hairer, BHZ, BCCH, CH, HS}, paracontrolled calculus \cite{GIP, BB1, BB2, BB3} or Duch's renormalization group setting \cite{Duch, DuchElliptic}. Either way one gets (in a Euclidean setting) from the use of any of these tools a proper definition of a solution to Equation \eqref{EqPhi43SPDE} and a local in time well-posedness result that needs to be supplemented by some ad hoc arguments to prove the long time existence of its solution. The Markovian character of the dynamics on $C^{-1/2-\epsilon}(M)$ generated Equation \eqref{EqPhi43SPDE} is inherited from its discrete counterpart. A compactness argument related to the property of `{\sl coming down from infinity}' satisfied by the solutions of Equation \eqref{EqPhi43SPDE} then gives the long-time existence of the local solution and the existence of an invariant measure for the semigroup on $C^{-1/2-\epsilon}(M)$ generated by this equation. This was first proved in the setting of the torus by Mourrat \& Weber in \cite{MourratWeber}. The uniqueness of such an invariant measure was proved in the $3$-dimensional torus using a robust argument from dynamical systems: If the semigroup generated by the dynamics \eqref{EqPhi43SPDE} has the strong Feller property and there is in the state space an accessible point then the semigroup has at most one invariant probability measure. Hairer \& Mattingly proved in \cite{HairerMattingly} a general result that shows in particular that the $\Phi^4_3$ dynamics on the $3$-dimensional torus has the strong Feller property. Hairer \& Sch\"onbauer proved in \cite{HairerSchonbauer} a very general and deep result on the support of the law of a certain class of random models that gives as a by-product the existence of an accessible point for the $\Phi^4_3$ dynamics on the $3$-dimensional torus. None of these results are available in a manifold setting, and proving them in a manifold setting in the same generality as in \cite{HairerMattingly} and \cite{HairerSchonbauer} appears to us as a considerable task.

\ssk

We did not use regularity structure, paracontrolled calculus or renormalization group methods in our construction of `a' $\Phi^4_3$ measure on an arbitrary boundaryless $3$-dimensional Riemannian manifold $M$ in \cite{BDFT}. Rather we followed Jagannath \& Perkowski who noticed in \cite{JagannathPerkowski} that a clever change of variable allows to rewrite the proper formulation of Equation \eqref{EqPhi43SPDE} as a PDE with random coefficients
\begin{equation} \label{EqJPFormulationSynthetic}
(\partial_t-\Delta+1)v = B\cdot\nabla v - Av^3 + Z_2v^2 + Z_1v + Z_0,
\end{equation}
where $B\in C([0,T], C^{-\eta}(M,TM)), A\in C([0,T], C^{1-\eta}(M))$ and
$$
Z_i\in C([0,T], C^{-\frac{1}{2}-\eta}(M)) \qquad (0\leq i\leq 2)
$$ 
are random variables built from the noise, for any $0<T<\infty$ and $\eta>0$. The dynamics \eqref{EqJPFormulationSynthetic} is completed with the datum of an initial condition in a space of the form $C^{-1/2-\epsilon}(M)$, for $\epsilon>0$ small enough. Note that no singular product is involved in Equation \eqref{EqJPFormulationSynthetic}; the renormalization problem in \eqref{EqPhi43SPDE} is involved in the definition and construction of the random variables $A, Z_2,Z_1,Z_0$. 

We were able in \cite{BDFT} to construct these random fields and prove an $L^p$ coming down from infinity result for the solutions to Equation \eqref{EqJPFormulationSynthetic} that entails the existence of an invariant probability measure for the Markovian dynamics \eqref{EqPhi43SPDE}. The question of uniqueness of such an invariant probability measure was left aside in \cite{BDFT}; this is the point that we address in the present work.

\ssk

\begin{thm} \label{ThmMainInvariance}
The semigroup on $C^{-1/2-\epsilon}(M)$ generated by the dynamics \eqref{EqPhi43SPDE} has a unique invariant probability measure.
\end{thm}

\ssk

The change of variable $(u\mapsto v)$ from \eqref{EqPhi43SPDE} to \eqref{EqJPFormulationSynthetic} is explicit: Adding for instance a (possibly random adapted) drift $h$ in the dynamics of $u$ adds an explicit $h$-dependent drift in the dynamics of $v$. We use Equation \eqref{EqJPFormulationSynthetic} as a convenient description of the Markovian dynamics of $u$ to construct a coupling by change of measure between two solutions of Equation \eqref{EqPhi43SPDE} started from two arbitrary initial conditions in the state space. We obtain some explicit control on the probability of a successful coupling that is independent of the pair of initial conditions. This allows us to infer the uniqueness of an invariant probability measure for the dynamics generated by \eqref{EqPhi43SPDE}. To run this approach we need to strengthen the $L^p$ coming down from infinity result proved in \cite{BDFT} into an $L^\infty$ coming down result for the solution $v$ to Equation \eqref{EqJPFormulationSynthetic}. This is what Section \ref{SectionLinftyComingDown} is about. We adapt there to our setting Moinat \& Weber' seminal approach \cite{MoinatWeber} to the coming down phenomenon. This kind of control is actually needed not only for $v$ but also for the solution $v_\ell$ of an equation similar to Equation \eqref{EqPhi43SPDE}, with an additional drift that depends on a real parameter $\ell$. Section \ref{SectionVariationTheme} deals with that perturbed equation. As a matter of fact it turns out to be necessary to also have a quantitative control on the sizes of $v(t)$ and $v_\ell(t)$ in stronger norms, not just $L^\infty$; such controls are provided in Section \ref{SubsectionControlNorm32} and Section \ref{SectionVariationTheme}. Equipped with the quantitative estimates proved in these sections we construct in Section \ref{SectionUniqueness} a coupling by change of measure that leads to a proof of uniqueness of an invariant measure for the semigroup generated by \eqref{EqPhi43SPDE}. As a by-product of our analysis we prove in Section \ref{SubsectionStrongFeller} a Harnack-type inequality for the semigroup that provides a short proof that this semigroup has the strong Feller property. 

\ssk

Our uniqueness result gives a characterization of our $\Phi^4_3$ measure as the unique invariant probability measure of a Markovian dynamics on some distribution space over $M$. As this dynamics depends only on the Riemannian structure of $M$ the $\Phi^4_3$ measure appears as depending only on the isometry class of the Riemannian manifold $M$. As a matter of fact, the reasoning that we use to prove the uniqueness of the $\Phi^4_3$ measure is insensitive to the geometry of $M$; the latter plays no role in what follows. The $\Phi^4_3$ measure itself is sensitive to the metric and it would be interesting to see which features of the measure are influenced by the curvature of the metric for instance. Let us mention here that there is not a unique way of renormalizing Equation \eqref{EqPhi43SPDE}. There is actually a one parameter family of renormalized equations. We work here with a particular value of that parameter that was fixed in \cite{BDFT} by the renormalization process adopted therein to construct $A, Z_2, Z_1, Z_0$. The reader should thus keep in mind that we talk below of `the' $\Phi^4_3$ measure as the unique invariant probability measure of our choice of formulation of the dynamics \eqref{EqPhi43SPDE}, even though there is a one parameter family of $\Phi^4_3$ measure corresponding to a one parameter family of versions of Equation \eqref{EqPhi43SPDE}.

\ssk

Before this work there were two ways of proving the uniqueness of the $\Phi^4_3$ measure in the torus. 
\begin{itemize}
	\item[--] The strong Feller plus accessibility strategy described above.

	\item[--] Barashkov \& Gubinelli's stochastic control representation of the partition function \cite{BG}, based on the Bou\'e-Dupuis formula.
\end{itemize}
The ingredients \cite{HairerMattingly, HairerSchonbauer} of the first option are not available in a manifold setting. The construction of the stochastic objects involved in \cite{BG} was not achieved yet in a manifold setting and this is a non-trivial task. The coupling strategy that we implement in the present work does not require any of the arguments used in the two other strategies. It provides a third alternative in the flat torus case, where it also gives a new Harnack-type inequality. Couplings by change of measures have been largely used in the study of stochastic differential equations and It\^o type stochastic partial differential equations, especially by F.Y. Wang and his co-authours. They provide a robust path to a whole family of functional inequalities. We highly recommend his little book \cite{WangBook} on the subject.

\medskip

\noindent \textbf{Notation} -- {\it For an initial condition $\phi$ of \eqref{EqPhi43SPDE} we will denote by $\phi'$ the corresponding initial condition of \eqref{EqJPFormulationSynthetic} given by the Jagannath \& Perkowski transform
\begin{equation} \label{EqUFromV}
\phi = \begin{tikzpicture}[scale=0.3,baseline=0cm]
\node at (0,0)  [dot] (1) {};
\node at (0,0.8)  [noise] (2) {};
\draw[K] (1) to (2);
\end{tikzpicture}(0) - \begin{tikzpicture}[scale=0.3,baseline=0cm]
\node at (0,0) [dot] (0) {};
\node at (0,0.5) [dot] (1) {};
\node at (-0.4,1)  [noise] (noise1) {};
\node at (0,1.2)  [noise] (noise2) {};
\node at (0.4,1)  [noise] (noise3) {};
\draw[K] (0) to (1);
\draw[K] (1) to (noise1);
\draw[K] (1) to (noise2);
\draw[K] (1) to (noise3);
\end{tikzpicture}(0)
+ e^{-3\begin{tikzpicture}[scale=0.3,baseline=0cm]
\node at (0,0) [dot] (0) {};
\node at (0,0.4) [dot] (1) {};
\node at (-0.3,0.8)  [noise] (noise1) {};
\node at (0.3,0.8)  [noise] (noise2) {};
\draw[K] (0) to (1);
\draw[K] (1) to (noise1);
\draw[K] (1) to (noise2);
\end{tikzpicture}(0)} \big(\phi' + v_{\textrm{ref}}(0)\big),
\end{equation}
with the notations of \cite{BDFT}. Roughly speaking, $\begin{tikzpicture}[scale=0.3,baseline=0cm]
\node at (0,0)  [dot] (1) {};
\node at (0,0.8)  [noise] (2) {};
\draw[K] (1) to (2);
\end{tikzpicture}(0) = \int_{-\infty}^0 e^{s\Delta}(\xi) ds$, the random variables  $\begin{tikzpicture}[scale=0.3,baseline=0.05cm]
\node at (0,0) [dot] (1) {};
\node at (-0.3,0.5)  [noise] (noise1) {};
\node at (0.3,0.5)  [noise] (noise2) {};
\draw[K] (1) to (noise1);
\draw[K] (1) to (noise2);
\end{tikzpicture}$ and $\begin{tikzpicture}[scale=0.3,baseline=0.05cm]
\node at (0,0) [dot] (1) {};
\node at (-0.4,0.5)  [noise] (noise1) {};
\node at (0,0.7)  [noise] (noise2) {};
\node at (0.4,0.5)  [noise] (noise3) {};
\draw[K] (0) to (1);
\draw[K] (1) to (noise1);
\draw[K] (1) to (noise2);
\draw[K] (1) to (noise3);
\end{tikzpicture}$ are the Wick square and power $3$ of $\begin{tikzpicture}[scale=0.3,baseline=0.05cm]
\node at (0,0)  [dot] (1) {};
\node at (0,0.8)  [noise] (2) {};
\draw[K] (1) to (2);
\end{tikzpicture}$ and $\begin{tikzpicture}[scale=0.3,baseline=0.05cm]
\node at (0,0) [dot] (1) {};
\node at (0,1) [dot] (2) {};
\draw[K] (1) to (2);
\end{tikzpicture}(0)$ stands for the integration operator $\int_{-\infty}^0 e^{s\Delta}(\cdot) ds$. The precise definition of the different terms above plays no role here, so we refer the interested reader to Section 2 of \cite{JagannathPerkowski}.  In the reformulation \eqref{EqJPFormulationSynthetic} of the parabolic $\Phi^4_3$ equation \eqref{EqPhi43SPDE} the random variables $A,B,Z_2,Z_1,Z_0$ are constructed as limits in some appropriate spaces, as $r>0$ goes to $0$, of some random variables $A_r,B_r,Z_{2r},Z_{1r},Z_{0r}$ that are some explicit functionals of a regularized noise $\xi_r$. With the previous pictorial notation one has for instance $A_r = \exp\big(-6\begin{tikzpicture}[scale=0.3,baseline=0cm]
\node at (0,0) [dot] (0) {};
\node at (0,0.4) [dot] (1) {};
\node at (-0.3,0.8)  [noise] (noise1) {};
\node at (0.3,0.8)  [noise] (noise2) {};
\draw[K] (0) to (1);
\draw[K] (1) to (noise1);
\draw[K] (1) to (noise2);
\end{tikzpicture}_r\big)\in C([0,T],C^{1-\eta}(M))$ and $B=-6\nabla \begin{tikzpicture}[scale=0.3,baseline=0cm]
\node at (0,0) [dot] (0) {};
\node at (0,0.4) [dot] (1) {};
\node at (-0.3,0.8)  [noise] (noise1) {};
\node at (0.3,0.8)  [noise] (noise2) {};
\draw[K] (0) to (1);
\draw[K] (1) to (noise1);
\draw[K] (1) to (noise2);
\end{tikzpicture} \in  C([0,T],C^{-\eta}(M,TM))$. One can find in Equation (2.4) of Jagannath \& Perkowski's work \cite{JagannathPerkowski} the explicit, non-informing, expression of $Z_{2r},Z_{1r},Z_{0r}$. We note that the random variables $A,B,Z_2,Z_1,Z_0$ are adapted to the natural time filtration of $\xi$ and stationary in time. We use the notation $\widehat\xi$ for the tuple $(A,Z_2,Z_1,Z_0)$.   }

\medskip

\noindent \textbf{Acknowledgements} -- {\it The content of Section \ref{SectionLinftyComingDown} and Section \ref{SubsectionControlNorm32} owes a lot to discussions with T.D. T\^o. I also warmly thank V.N. Dang for his indirect input in this work.}

\bigskip

\section{Uniqueness of the invariant measure}
\label{SectionUniqueness}

The dynamics on $C^{-1/2-\epsilon}(M)$ generated by \eqref{EqPhi43SPDE} is Markovian and defined for all times. The global in time well-posedness comes from the strong damping effect of the $-u^3$ term, which plays a crucial role in obtaining some explicit control on the $L^p$ norm of the solution $v$ to Equation \eqref{EqJPFormulationSynthetic} that is independent of its initial condition. This control on $v$ then gives a control on $u$ that is independent of its initial condition, in an appropriate space, via the explicit correspondence $v\mapsto u$. This phenomenon is called `{\sl coming down from infinity}'. It was first proved for a transform of the solutions of Equation \eqref{EqPhi43SPDE} by Mourrat \& Weber in their seminal work \cite{MourratWeber} on the $\Phi^4$ equation on the $3$-dimensional torus. It was later extended to the Euclidean setting of $\bbR^3$ by Moinat \& Weber \cite{MoinatWeber} and Gubinelli \& Hofmanov\'a \cite{GubinelliHofmanova} using different methods. 

The coming down from infinity estimate can be used to get the existence of an invariant probability measure by an elementary compactness argument involving the time stationarity of some random fields. (This is for instance done in Section 5.1 of \cite{BDFT}.) We prove in this section that the Markovian dynamics \eqref{EqPhi43SPDE} has at most one invariant probability measure. We work for that purpose with Jagannath \& Perkowski's formulation \eqref{EqJPFormulationSynthetic} of \eqref{EqPhi43SPDE} and use a {\it coupling argument} to prove the uniqueness. Given two points $\phi'_1,\phi'_2\in C^{-1/2-\epsilon}(M)$, an elementary coupling of two solutions to \eqref{EqJPFormulationSynthetic} started from $\phi'_1$ and $\phi'_2$ would consist in constructing on some probability space a pair of spacetime white noises such that the solutions to Equation \eqref{EqJPFormulationSynthetic} built from each of these noises take the same value at some fixed positive finite time $T$ outside of an event of arbitrarily small probability independent of $\phi'_1,\phi'_2$. The corresponding solutions of Equation \eqref{EqPhi43SPDE} would also coincide at that time. One could then take some random initial conditions with law two invariant probability measures $\mu_1,\mu_2$ for the dynamics generated by \eqref{EqPhi43SPDE} and write for any bounded continuous function $f$ on $C^{-1/2-\epsilon}(M)$
$$
\mu_1(f) = \bbE\big[f(u(T ; \mu_1))\big] = \bbE\big[f(u(T ; \mu_2))\big] = \mu_2(f),
$$
with $u(\cdot\,;\,\mu_i)$ denoting the solution to  \eqref{EqPhi43SPDE} with random initial condition with law $\mu_i$. We are not able to produce such a strong coupling here; rather, given the trajectory $u(\cdot\,;\,\mu_1)$, we can add to the dynamics of $u(\cdot\,;\,\mu_2)$ a drift that forces the latter to meet the former by a fixed time $T$ with high probability. With the dynamics of $u(\cdot\,;\,\mu_2)$ changed the measure $\mu_2$ is not invariant anymore for this new dynamics and the above simple argument for uniqueness does not apply per se. However, for a particular drift there is an equivalent probability measure on our probability space for which the new dynamics has the same law as the original dynamics with random initial condition with law $\mu_2$. A variation on the above pattern of proof then gives the equality of $\mu_1$ and $\mu_2$. (This variation is detailed below in Step 2 of the proof of Theorem \ref{ThmUniqueness}.)

\medskip

Denote by $(\Omega, \mcF, \bbP)$ the probability space on which all our random variables have been implicitly defined so far. We write $\mathscr{L}_\bbP(X)$ for the law under $\bbP$ of a random variable $X$ and use a similar notation $\mathscr{L}_\bbQ(X)$ for any other probability measure $\bbQ$ on $(\Omega,\mcF)$. We write $\bbE_\bbP$ and $\bbE_\bbQ$ for the corresponding expectation operators. We use a coupling by a change of measure argument to prove that the semigroup $(\mcP_t)_{t\geq 0}$ on $C^{-1/2-\epsilon}(M)$ generated by \eqref{EqPhi43SPDE} has at most one invariant probability measure.

\ssk

\begin{thm} \label{ThmUniqueness}
The semigroup $(\mcP_t)_{t\geq 0}$ has at most one invariant probability measure.
\end{thm}

\ssk

We will need along the proof of Theorem \ref{ThmUniqueness} some quantitative results that will be proved later in Section \ref{SectionLinftyComingDown} and Section \ref{SubsectionControlNorm32}. These estimates allow us to use the result of Lemma \ref{LemMonotony} below.

\ssk

\begin{Dem}
We proceed in two steps and first construct a coupling by a change of measure between two trajectories of the Jagannath \& Perkowski version of the $\Phi^4_3$ dynamics started from different points.  As a preliminary remark note that shifting the noise $\xi$ by a (possibly random) element $h$ of its Cameron-Martin space with support in time in the interval $[1,2]$ is equivalent to adding a drift $h \exp\big(3\begin{tikzpicture}[scale=0.3,baseline=0cm] \node at (0,0) [dot] (0) {}; \node at (0,0.4) [dot] (1) {}; \node at (-0.3,0.8)  [noise] (noise1) {}; \node at (0.3,0.8)  [noise] (noise2) {}; \draw[K] (0) to (1); \draw[K] (1) to (noise1); \draw[K] (1) to (noise2); \end{tikzpicture}\big)$ to the dynamics of $v$. Indeed, let $\phi$ and $\phi'$ be related by the relation 
$$
\phi' = \exp\big(3\begin{tikzpicture}[scale=0.3,baseline=0cm] \node at (0,0) [dot] (0) {}; \node at (0,0.4) [dot] (1) {}; \node at (-0.3,0.8)  [noise] (noise1) {}; \node at (0.3,0.8)  [noise] (noise2) {}; \draw[K] (0) to (1); \draw[K] (1) to (noise1); \draw[K] (1) to (noise2); \end{tikzpicture}(0)\big) \Big(\phi  -  \begin{tikzpicture}[scale=0.3,baseline=0cm] \node at (0,0)  [dot] (1) {}; \node at (0,0.8)  [noise] (2) {}; \draw[K] (1) to (2); \end{tikzpicture}(0) + \begin{tikzpicture}[scale=0.3,baseline=0cm] \node at (0,0) [dot] (0) {}; \node at (0,0.5) [dot] (1) {}; \node at (-0.4,1)  [noise] (noise1) {}; \node at (0,1.2)  [noise] (noise2) {}; \node at (0.4,1)  [noise] (noise3) {}; \draw[K] (0) to (1); \draw[K] (1) to (noise1); \draw[K] (1) to (noise2); \draw[K] (1) to (noise3); \end{tikzpicture}(0)\Big) - v_{\textrm{ref}}(0).
$$
One sees that if $v_r^h$ solves an equation of the form

$$
(\partial_t-\Delta)v_r^h = B_r\cdot\nabla v_r^h - A_r(v_r^h)^3 + Z_{2,r}(v_r^h)^2 + Z_{1,r} v_r^h + Z_{0,r} + (e^{r\Delta} h) \, e^{3\begin{tikzpicture}[scale=0.3,baseline=0cm] \node at (0,0) [dot] (0) {}; \node at (0,0.4) [dot] (1) {}; \node at (-0.3,0.8)  [noise] (noise1) {}; \node at (0.3,0.8)  [noise] (noise2) {}; \draw[K] (0) to (1); \draw[K] (1) to (noise1); \draw[K] (1) to (noise2); \end{tikzpicture}_r},
$$
with initial condition $v_r^h(0)=\phi'$, then 
\begin{equation} \label{EqJPTransform}
u_r^h = \begin{tikzpicture}[scale=0.3,baseline=0cm] \node at (0,0)  [dot] (1) {}; \node at (0,0.8)  [noise] (2) {}; \draw[K] (1) to (2); \end{tikzpicture}_r - \begin{tikzpicture}[scale=0.3,baseline=0cm] \node at (0,0) [dot] (0) {}; \node at (0,0.5) [dot] (1) {}; \node at (-0.4,1)  [noise] (noise1) {}; \node at (0,1.2)  [noise] (noise2) {}; \node at (0.4,1)  [noise] (noise3) {}; \draw[K] (0) to (1); \draw[K] (1) to (noise1); \draw[K] (1) to (noise2); \draw[K] (1) to (noise3); \end{tikzpicture}_r + e^{-3\begin{tikzpicture}[scale=0.3,baseline=0cm] \node at (0,0) [dot] (0) {}; \node at (0,0.4) [dot] (1) {}; \node at (-0.3,0.8)  [noise] (noise1) {}; \node at (0.3,0.8)  [noise] (noise2) {}; \draw[K] (0) to (1); \draw[K] (1) to (noise1); \draw[K] (1) to (noise2); \end{tikzpicture}_r} (v_r^h + v_{\textrm{ref},r}),
\end{equation}
solves the equation
$$
(\partial_t-\Delta)u_r^h = -(u_r^h)^3 + 3(a_r-b_r)u_r^h + e^{r\Delta}(\xi+h), \qquad u_r(0)=\phi.
$$
The convergence of $v_r^h$ to the solution $v^h$ of the equation
$$
(\partial_t-\Delta)v^h = B\cdot\nabla v^h - A(v^h)^3 + Z_2 (v^h)^2 + Z_1 v^h + Z_0 + h e^{3\begin{tikzpicture}[scale=0.3,baseline=0cm] \node at (0,0) [dot] (0) {}; \node at (0,0.4) [dot] (1) {}; \node at (-0.3,0.8)  [noise] (noise1) {}; \node at (0.3,0.8)  [noise] (noise2) {}; \draw[K] (0) to (1); \draw[K] (1) to (noise1); \draw[K] (1) to (noise2); \end{tikzpicture}}, \qquad v^h(0)=\phi',
$$
ensures the convergence of $u_r^h$ to a limit.

\ssk

\noindent $\rhd$ \textbf{\textsf{Step 1  -- The coupling.}} Pick $\phi_1,\phi_2\in C^{-1/2-\epsilon}(M)$ with corresponding $\phi_1',\phi_2'$. We adopt as above the notation $v = v(\cdot,\phi'_1)$ for the solution of the Jagannath-Perkowski equation with initial condition $\phi_1'$, with $u=u(\cdot,\phi_1)$ the corresponding function given by the inverse Jagannath-Perkowski transform. Let $v_\ell$ be the solution of the equation
\begin{equation} \label{EqVEll} \begin{split}
\partial_tv_\ell &= (\Delta-1)v_\ell - Av_\ell^3 + B\cdot\nabla v_\ell + Z_2v_\ell^2 + Z_1v_\ell + Z_0   \\
&\quad+ \ell \, {\bf 1}_{1< t <\tau}\,\frac{v(t)-v_\ell(t)}{\Vert v(t)-v_\ell(t)\Vert_{L^2}}\,\exp\big(3\begin{tikzpicture}[scale=0.3,baseline=0cm] \node at (0,0) [dot] (0) {}; \node at (0,0.4) [dot] (1) {}; \node at (-0.3,0.8)  [noise] (noise1) {}; \node at (0.3,0.8)  [noise] (noise2) {}; \draw[K] (0) to (1); \draw[K] (1) to (noise1); \draw[K] (1) to (noise2); \end{tikzpicture}(t)\big)   \\
&\eqdef F(v_\ell) + \ell \, {\bf 1}_{1< t <\tau}\,\frac{v(t)-v_\ell(t)}{\Vert v(t)-v_\ell(t)\Vert_{L^2}} \, \exp\big(3\begin{tikzpicture}[scale=0.3,baseline=0cm] \node at (0,0) [dot] (0) {}; \node at (0,0.4) [dot] (1) {}; \node at (-0.3,0.8)  [noise] (noise1) {}; \node at (0.3,0.8)  [noise] (noise2) {}; \draw[K] (0) to (1); \draw[K] (1) to (noise1); \draw[K] (1) to (noise2); \end{tikzpicture}(t)\big) \eqdef F_\ell(t,v,v_\ell), \quad (0\leq t < \tau),
\end{split} \end{equation}
with initial condition $\Phi'_2$, and where 
$$
\tau = \tau(\ell,\phi_1,\phi_2) \defeq \inf\big\{s\geq 1 \, ; \, v_\ell(s) = v(s)\big\} \wedge 2.
$$
The random time $\tau$ is called the {\sl coupling time} -- we take as a convention $\inf\emptyset = + \infty$. A {\sl successful coupling} corresponds to the event $\{\tau<2\}$, in which case we let $v_\ell(t)=v(t)$ for $t\geq \tau$. We have
 
$$
{\bf 1}_{\tau<2} \, v_\ell(2, \phi'_2) = {\bf 1}_{\tau<2} \, v(2, \phi'_1)
$$ 
and 
$$
{\bf 1}_{\tau<2} \, u_\ell(2, \phi_2) = {\bf 1}_{\tau<2} \, u(2, \phi_1),
$$
with $u_\ell$ corresponding to $v_\ell$ via the inverse Jagannath-Perkowski transform \eqref{EqJPTransform}. Theorem \ref{ThmVariationOnTheme} below guarantees the longtime well-posedness of Equation \eqref{EqVEll} and provides some controls of $v_\ell$ of the form
$$
\Vert v_\ell(s)\Vert_\infty \lesssim_{\widehat{\xi}} 1+\ell^{1/3}
$$
and
$$
\Vert v_\ell(s)\Vert_{C^{1+2\epsilon}} \lesssim_{\widehat{\xi}} 1+\ell
$$
for all $1\leq s<\tau$. These estimates ensure that the assumptions involved in the next statement are relevant for the study of the joint dynamics of $v$ and $v_\ell$. Recall $\widehat\xi$ stands for the tuple $(A,Z_2,Z_1,Z_0)$.

\ssk

\begin{lem} \label{LemMonotony}
Take $\ell\geq 1$. There is an absolute constant $\epsilon_0>0$ such that for all $0<\epsilon<\epsilon_0$ and for $1<t<\tau$ one has
\begin{equation} \label{EqUpperBoundMonotony}
\big\langle F(v) - F_\ell(t,v,v') , v - v'\big\rangle_{L^2} \lesssim_{\widehat{\xi}} - \big(\ell {\color{black} -} o(\widehat\xi ; \ell)\big) \Vert v - v'\Vert_{L^2}
\end{equation}
for all $v,v'\in C^{\color{black} 1+2\epsilon}(M)$ with 
$$
\Vert v\Vert_{L^\infty} \vee \Vert v'\Vert_{L^\infty} \leq c_1(\widehat\xi\,) \, {\color{black} \ell^{\frac{1}{3}}},
$$
and
\begin{equation} \label{Eq32SizeCondition}
{\color{black} \Vert v\Vert_{C^{1+2\epsilon}}  \vee \Vert v'\Vert_{C^{1+2\epsilon}} } \leq c_2(\widehat\xi\,) \, {\color{black} \ell},
\end{equation}
for a $\widehat\xi$-dependent {\color{black} non-negative} function $o(\widehat\xi ; \ell)$ of $\ell$ such that $o(\widehat\xi ; \ell)/\ell$ goes to $0$ as $\ell$ goes to $\infty$.
\end{lem}

\ssk

\textbf{\textsf{Proof --}} Note that there are no absolute values in \eqref{EqUpperBoundMonotony}. We give some upper bounds for each term in this expression. We make the common abuse of notation of writing $\langle Z_i f, g \rangle_{L^2}$ for $\langle Z_i fg , {\bf 1}\rangle$, the result of testing a well-defined distribution $Z_ifg$ on the constant function $\bf 1$. We use a similar convention for $\langle B\cdot\nabla f , g\rangle_{L^2}$.

\ssk

{\it -- The $A$ term.} As $A$ is positive one has
$$
\big\langle -A\big(v^3-{v'}^3\big) , v-v' \big\rangle_{L^2} \leq 0,
$$
and this term does not contribute to the upper bound \eqref{EqUpperBoundMonotony}. 

\ssk

{\it -- The $B$ term.} {\color{black} We start from the identity
$$
(v - v')^2 = 2(v-v') \prec (v-v') + (v-v') \odot (v-v'),
$$
with the left $(v-v')$ seen as an element of $L^2(M)$ and the right $(v-v')$ seen as an element of {\color{black} $C^{1+2\epsilon}(M)$} in the paraproduct and resonant terms, and estimate each term in $B^{1+\epsilon}_{2,\infty}(M)$. Losing a little bit on the regularity exponent allows using the interpolation size estimate between different Besov spaces and estimate
$$
\Vert v-v'\Vert_{B^{1+\epsilon}_{\infty\infty}} \lesssim_{\widehat\xi} \ell^{\frac{1}{3}\frac{\epsilon}{1+2\epsilon} + \frac{1+\epsilon}{1+2\epsilon}},
$$
with an exponent strictly smaller than $1$. We write
$$
\Vert v-v'\Vert_{B^{1+\epsilon}_{\infty\infty}} \lesssim_{\widehat\xi} \ell^{<1}.
$$
One then gets from the classical continuity estimates on the paraproduct and resonant operators that
$$
\Vert (v-v')^2\Vert_{B^{1+\epsilon}_{2\infty}} \lesssim_{\widehat\xi} \ell^{<1} \, \Vert v-v'\Vert_{L^2}.
$$

\ssk

{\it -- The $Z_1$ term.} We take advantage of the fact that $Z_1\in C_TB^{-\frac{1}{2}-\frac{\epsilon}{2}}_{21}(M)$ almost surely. We use the previous estimate to see that
$$
\big\vert \langle Z_1(v-v') , v-v'\rangle \big\vert \lesssim_{\widehat\xi, Z_1} \ell^{<1} \, \Vert v-v'\Vert_{L^2}.
$$

\ssk

{\it -- The $Z_2$ term.} Here as well we consider $Z_2$ as an element of $C_TB^{-\frac{1}{2}-\frac{\epsilon}{2}}_{21}(M)$. First we obtain by interpolation between the $L^\infty$ and $C^{1+2\epsilon}$ estimate on $v$ and $v'$ that
$$
\Vert v\pm v'\Vert_{B^{\frac{1}{2}+2\epsilon}_{\infty\infty}} \lesssim_{\widehat\xi} \ell^{\frac{1}{3}\frac{5+6\epsilon}{3-2\epsilon}},
$$
with an exponent slightly bigger than $5/9$. We thus get from the classical continuity estimates on the paraproduct and resonant operators that
\begin{equation} \label{EqSizeVMoinsVPrimeSquare}
\Vert (v-v')^2\Vert_{B^{\frac{1}{2}+2\epsilon}_{2\infty}} \lesssim_{\widehat\xi} \ell^{\frac{1}{3}\frac{5+6\epsilon}{3-2\epsilon}} \, \Vert v-v'\Vert_{L^2}.
\end{equation}
Now write
\begin{equation*} \begin{split}
\big(v^2-&(v')^2\big)(v-v')   \\
&= (v-v')^2(v+v')   \\
&= (v-v')^2 \prec (v+v') + \Big\{(v+v')\prec (v-v')^2 + (v+v')\odot (v-v')^2\Big\}.
\end{split} \end{equation*}
To estimate the contribution of the first paraproduct in the $Z_2$ term we use the elementary refined continuity estimate from Lemma \ref{LemSmallFactorTrick} in Appendix {\sf \ref{SectionAppendix}} to get the best of the $L^\infty$ and $C^{1+2\epsilon}$ estimates on $(v+v')$. We have for all integers $N$
\begin{equation*} \begin{split}
\big\Vert (v-v')^2 &\prec (v+v')\big\Vert_{B^{(\frac{1}{2}+2\epsilon)-\epsilon}_{2\infty}}   \\
&\lesssim_{\widehat\xi} \Vert (v-v')^2\Vert_{B^{\frac{1}{2}+2\epsilon}_{2\infty}}\Big( 2^{-N\epsilon} \Vert v+v'\Vert_{B^{\frac{1}{2}+2\epsilon}_{\infty\infty}} + N \Vert v+v'\Vert_{L^\infty}\Big)   \\
&\lesssim_{\widehat\xi} \ell^{\frac{1}{3} \frac{5+6\epsilon}{3-2\epsilon}} \, \Vert v-v'\Vert_{L^2} \Big( 2^{-N\epsilon} \ell^{\frac{1}{3} \frac{5+6\epsilon}{3-2\epsilon}} + N \ell^{\frac{1}{3}}\Big).
\end{split} \end{equation*}
Choosing $N$ such that $2^{-N\epsilon} \ell^{\frac{1}{3}\frac{5+6\epsilon}{3-2\epsilon}} \simeq \ell^{\frac{1}{3}}$ gives 
$$
2^{-N\epsilon/2} \ell^{\frac{1}{3}\frac{5+6\epsilon}{3-2\epsilon}} + N \ell^{\frac{1}{3}} \lesssim \ell^{\frac{1}{3}+\eta}
$$
for every $\eta>0$ and $\ell\geq \ell(\eta)$ large enough. One thus has
$$
\big\Vert (v-v')^2 \prec (v+v')\big\Vert_{B^{\frac{1}{2}+\epsilon}_{2\infty}} \lesssim_{\widehat\xi} \ell^{<1} \, \Vert v-v'\Vert_{L^2} 
$$
for an exponent $\frac{1}{3}\frac{5+6\epsilon}{3-2\epsilon} + \frac{1}{3} + \eta$ of $\ell$ strictly smaller than $1$, for an appropriate choice of $\eta$.

We can directly use \eqref{EqSizeVMoinsVPrimeSquare} and the $L^\infty$ estimate on $v$ and $v'$ to see that
$$
\big\Vert (v+v')\prec (v-v')^2 + (v+v')\odot (v-v')^2\big\Vert_{B^{\frac{1}{2}+2\epsilon}_{2\infty}} \lesssim_{\widehat\xi} \ell^{\frac{1}{3}\frac{5+6\epsilon}{3-2\epsilon} + \frac{1}{3}} \, \Vert v-v'\Vert_{L^2},
$$ 
here again with an exponent of $\ell$ strictly smaller than $1$.   }

\ssk

{\it -- The $Z_0$ term.} As we have a term $-\ell \Vert v-v'\Vert_{L^2}$ that comes from $F_\ell(t,v,v')$ and all the other contributions to $\big\langle F(v) - F_\ell(t,v,v') , v - v'\big\rangle_{L^2}$ add up to a quantity bounded above by a constant multiple of $\ell^{<1} \, \Vert v-v'\Vert_{L^2}$ we obtain \eqref{EqUpperBoundMonotony} for an appropriate choice of $\epsilon_0>0$.   \hspace{0.1cm}\hfill $\rhd$

\medskip

The proof makes it clear that one can take $o(\widehat\xi ; \ell)$ of the form
\begin{equation} \label{EqFormulaPetitOEll} \begin{split}
o(\widehat\xi ; \ell) = c\Big(\Vert B_{\vert[1,2]}\Vert_{C([1,2],C^{-\epsilon}(M))} &+ \Vert Z_{2\vert[1,2]}\Vert_{C([1,2],C^{-1/2-\epsilon})(M)}   \\
&+ \Vert Z_{1\vert[1,2]}\Vert_{C([1,2],C^{-1/2-\epsilon})(M)}\Big) \ell^\gamma
\end{split} \end{equation}
for some positive constant $c$ and some positive exponent $\gamma<1$. We will denote below by $m(\widehat\xi\,)$ the implicit multiplicative constant in $\lesssim_{\widehat\xi}$ in \eqref{EqUpperBoundMonotony}. It only depends on the restriction of $\widehat\xi$ to the time interval $[0,2]$. As in Lemma 4 of \cite{BDFT} one proves that the function of time $\Vert v(\cdot, \phi_1') - v_\ell(\cdot, \phi_2')\Vert_{L^2}$ is Young differentiable on the interval $(1,\tau)$ and one has for $1<t<\tau$
\begin{equation*} \begin{split}
\big\Vert v(t, \phi_1') &- v_\ell(t, \phi_2')\big\Vert_{L^2} - \big\Vert v(1, \phi_1') - v_\ell(1, \phi_2')\big\Vert_{L^2}   \\
&= \int_1^t \frac{\Big\langle F\big(v(s, \phi_1')\big) - F_\ell\big(t,v(s, \phi_1'),v_\ell(s, \phi_2')\big) , v(s, \phi_1') - v_\ell(s, \phi_2')\Big\rangle_{L^2}}{\big\Vert v(s, \phi_1') - v_\ell(s, \phi_2')\big\Vert_{L^2}}\,ds.
\end{split} \end{equation*}
It follows from Theorem \ref{ThmVariationOnTheme} and Lemma \ref{LemMonotony} that one has for $1\leq t <\tau$ the inequality
\begin{equation} \label{EqMajoration}
\big\Vert v(t, \phi_1') - v_\ell(t, \phi_2')\big\Vert_{L^2} \leq \big\Vert v(1,\phi'_1) - v_\ell(1,\phi'_2)\big\Vert_{L^2} - m(\widehat\xi\,)\big(\ell - o(\widehat\xi ; \ell)\big) (t-1),
\end{equation}
for a positive quantity $o(\widehat\xi ; \ell)$ that depends on $\widehat\xi_{\vert [1,2]}$ such that $o(\widehat\xi ; \ell)/\ell$ goes to $0$ as $\ell$ goes to $\infty$. (Theorem \ref{ThmVariationOnTheme} ensures that one can use Lemma \ref{LemMonotony}.) So one has a successful coupling on the event
$$
\bigg\{\big\Vert v(1,\phi'_1) - v_\ell(1,\phi'_2) \big\Vert_{L^2} \leq m(\widehat\xi\,)\,\frac{\ell - o(\widehat\xi ; \ell)}{2} \bigg\} \subset \big\{\tau<2\big\}.
$$
As a consequence of this inclusion and the $L^p$ or $L^\infty$ coming down from infinity result of \cite{BDFT} or Theorem \ref{ThmLinftyComingDown}, one can choose $\ell$ big enough to have both $\bbP(\tau(\ell, \phi_1, \phi_2)=2)$ and $\bbQ_\ell(\tau(\ell,\phi_1, \phi_2)=2)$ strictly smaller than $1$ independently of $\phi_1,\phi_2$, say
\begin{equation} \label{EqUniformBoundNonCoupling}
\max\big(\bbP(\tau(\ell,\phi_1', \phi_2')=2), \bbQ_\ell(\tau(\ell,\phi_1', \phi_2')=2)\big)\leq a < 1,
\end{equation} 
for all $\phi'_1,\phi'_2\in C^{-1/2-\epsilon}(M)$. We fix such an $\ell$ and set
\begin{equation*}
R_{\ell,\phi_1,\phi_2} \defeq \exp\bigg(- \ell \, \xi\Big({\bf 1}_{1<\cdot<\tau} \,\frac{v(\cdot,\phi'_1) - v_\ell(\cdot,\phi'_2)}{\Vert v(\cdot,\phi'_1) - v_\ell(\cdot,\phi'_2) \Vert_{L^2}}\Big) - \frac{\ell^2{\color{black} (\tau-1)}}{2}\bigg).
\end{equation*}
Note that $v$ and $v_\ell$ are both adapted to the filtration of $\xi$. Since $\tau\leq 2$, Novikov's integrability criterion
\begin{equation*}
\bbE\Big[\exp\Big(\frac{\ell^2{\color{black} (\tau-1)}}{2}\Big)\Big] < \infty
\end{equation*}
is satisfied and it follows from Girsanov theorem that the process
$$
\xi + \ell \, {\bf 1}_{1<\cdot<\tau}\,\frac{v - v_\ell}{\Vert v - v_\ell\Vert_{L^2}}
$$
is under the probability
$$
d\bbQ_{\ell,\phi_1,\phi_2} \defeq R_{\ell,\phi_1,\phi_2} \, d\bbP
$$
a spacetime white noise. (See e.g. Theorem 10.14 of Da Prato \& Zabczyk \cite{DaPratoZabczyk} for that type of statement.) Pick $\alpha\in(0,1]$. We have 
\begin{equation} \label{EqIdentityInLaw}
\mathscr{L}_{\bbQ_{\ell,\phi_1,\phi_2}}(u_\ell(\cdot,\phi_2)) = \mathscr{L}_{\bbP}(u(\cdot,\phi_2))
\end{equation}
and 
\begin{equation} \label{EqCouplingIdentity}
u_\ell(2, \phi_2) = u(2,\phi_1) \textrm{ on the event } \big\{\tau < 2\big\}.
\end{equation}

\ssk

\noindent $\rhd$ \textbf{\textsf{Step 2 -- Uniqueness of an invariant probability measure.}} We can now prove that the semigroup $(\mcP_t)_{t\geq 0}$ has at most one invariant probability measure. Otherwise, there would be (at least) two extremal invariant, hence singular, probability measures $\mu,\nu$. (See e.g. Theorem 1.7 in Hairer's lecture notes \cite{LNHairerMarkov} on the convergence of Markov processes.) We could take $\phi_1$ random with law $\mu$, and $\phi_2$ random with law $\nu$, and keep writing $\bbE$ for the expectation operator in this extended probability space. Simply write $R_{\ell}$ rather than $R_{\ell,\phi_1,\phi_2}$. Write $u_{\ell}(\cdot,\nu)$ and $u(\cdot,\mu)$ to emphasize the law of the initial condition. {\color{black} For a measurable set $A\subset C^{-1/2-\epsilon}(M)$ with $\mu(A) = 0$ we prove below that $\nu(A)=0$. The measure $\nu$ would thus be absolutely continuous with respect to $\mu$, a contradiction with the fact that $\nu$ is singular with respect to $\mu$. We give two proofs.

\ssk

{\it 1.} Assuming $\mu(A) = 0$, one would have from the identity in law \eqref{EqIdentityInLaw} and the fact that 
$$
\bbP\big(u(2,\mu) \in A\big) = \mu(\mcP_t{\bf 1}_A) = \mu(A) = 0
$$ 
the identity
\begin{equation*} \begin{split}
\nu(A) &= \nu(\mcP_2{\bf 1}_A) = \bbQ_{\ell}\big(u_{\ell}(2,\nu)\in A\big) = \bbE\big[R_{\ell} {\bf 1}_A(u_{\ell}(2,\nu))\big]   \\
&\hspace{-0.15cm}\overset{\eqref{EqCouplingIdentity}}{=} \bbE\big[R_{\ell} {\bf 1}_A(u(2,\mu)) {\bf 1}_{\tau<2}\big] + \bbE\big[R_{\ell} {\bf 1}_A(u_{\ell}(2,\nu)) {\bf 1}_{\tau=2}\big]    \\
&= \bbE\big[R_{\ell} {\bf 1}_A(u_{\ell}(2,\nu)) {\bf 1}_{\tau=2}\big] \leq \bbQ_{\ell}(\tau=2) \overset{\eqref{EqUniformBoundNonCoupling}}{\leq} a < 1.
\end{split} \end{equation*}
This is not enough to conclude that $\nu(A)=0$, but instead of coupling the two dynamics on a single time interval $[1,2]$ we can repeat if necessary our attempts to couple them a fixed finite number of times, during the time intervals $[2k-1,2k]$ after a coupling-free evolution on the time interval $[2k-2,2k-1]$, for $k\leq n$, say. Denote by $u_\ell^{(n)}(\cdot,\phi_2)$ the corresponding dynamics. Write $\tau_1\in [1,2], \dots, \tau_n\in[n+1,n+2]$ for the successive coupling times and set 
$$
\ln R^{(n)}_{\ell, \phi_1,\phi_2} \defeq -\hspace{-0.07cm}\sum_{k=1}^n \bigg(\ell \, \xi\Big({\bf 1}_{2k-1<\cdot<\tau_k} \,\frac{v(\cdot,\phi'_1) - v_\ell(\cdot,\phi'_2)}{\Vert v(\cdot,\phi'_1) - v_\ell(\cdot,\phi'_2) \Vert_{L^2}}\Big) + \frac{\ell^2 \big(\tau_k-2k+1\big)}{2}\bigg)
$$
and 
$$
d\bbQ_\ell^{(n)} \defeq R_\ell^{(n)}d\bbP.
$$
The probability measure $\bbQ_\ell^{(n)}$ implicitly depends on $\phi_1$ and $\phi_2$ and we also set
$$
\overline\bbQ_\ell^{(n)} \defeq \int \bbQ_\ell^{(n)}\nu(d\phi_2)\mu(d\phi_1).
$$
The process $u_\ell^{(n)}(\cdot,\phi_2)$ has under $\bbQ_\ell^{(n)}$ the same distribution as $u(\cdot,\phi_2)$ and the pair 
$$
\big(u(\cdot,\phi_1), u_\ell^{(n)}(\cdot,\phi_2)\big)
$$ 
is Markovian under both $\bbP$ and $\bbQ_\ell^{(n)}$. Denote by $\theta_r : \Omega\rightarrow\Omega$ a family of measurable measure preserving maps on $(\Omega,\mcF, \bbP)$ such that $\theta_r\circ\theta_{r'}=\theta_{r+r'}$ and $\theta_r \xi(\cdot) = \xi(\cdot+r)$. This shift acts on measurables functions of $\xi$ such as $u$ and $u_\ell^{(n)}$. One then has as above

\begin{equation*} \begin{split}
\nu(A) &= \int \bbE_\ell^{(n)}\Big[{\bf 1}_A\big(u_\ell^{(n)}(2n, \phi_2)\big){\bf 1}_{\tau_n(\phi_1,\phi_2)=2n+2}\Big]\nu(d\phi_2)\mu(d\phi_1)   \\
	  &= \int \bbE_\ell^{(n)}\Big[{\bf 1}_A\big(\theta_{2n-2}u_\ell(2, u_\ell^{(n-1)}(2n-2,\phi_2))\big){\bf 1}_{\tau_n(\phi_1,\phi_2)=2n+2}\Big]\nu(d\phi_2)\mu(d\phi_1)   \\
	  &= \int \bbE_\ell^{(n-1)}\Big[ \bbE_\ell^{(1)}\Big[ {\bf 1}_A\big(\theta_{2n-2}u_\ell^{(1)}(2, u_\ell^{(n-1)}(2n-2,\phi_2))\big) {\bf 1}_{\tau_1=2} \Big\vert u_\ell^{(n-1)}(2n-2,\phi_2) \Big]   \\ 
	  &\hspace{6cm} \times{\bf 1}_{\tau_{n-1}(\phi_1,\phi_2)=2n}\Big]\nu(d\phi_2)\mu(d\phi_1)   \\
	  &\leq a\,\overline\bbQ_\ell^{(n-1)}(\tau_{n-1}=2) \leq a^n,
\end{split} \end{equation*}
by induction. The conclusion $\nu(A)=0$ follows from the fact that $n$ is arbitrary.  

\ssk

{\it 2.} Alternatively one can assume that the two invariant probability measures $\mu$ and $\nu$ are singular and proceed as follows to get a contradiction. Denote by $\bbP_1$ the law of $u(\cdot,\mu)$ and by $\bbP_2$ the law of $u(\cdot,\nu)$, with a time parameter running in the time interval $[0,2]$. Step 1 produces a coupling between $\bbP_1$ and a probability with positive density $D$ with respect to $\bbP_2$. This coupling is a probability measure $\bbQ$ on $C\big([0,2],C^{-1/2-\epsilon}(M)\big)\times C\big([0,2],C^{-1/2-\epsilon}(M)\big)$ that gives a positive probability to the event $\big\{u_1(2)=u_2(2)\big\}$, denoting by $u_1$ and $u_2$ the canonical marginal processes on the product space. Denote by $\pi_1$ and $\pi_2$ the canonical projections and set
\begin{equation} \label{EqQPlus}
d\bbQ^- \defeq (1\wedge D^{-1})d\bbQ
\end{equation}
and 
$$
\bbQ^+ \defeq \bbQ^- + (\bbP_1-\pi_{1\star}\bbQ^-)\otimes (\bbP_2-\pi_{2\star}\bbQ^-).
$$
The measure $\bbQ^+$ on $C\big([0,2],C^{-1/2-\epsilon}(M)\big)\times C\big([0,2],C^{-1/2-\epsilon}(M)\big)$ is a probability measure with marginals $\bbP_1$ and $\bbP_2$ -- that is, a coupling of these two probability measures. We further have that $\bbQ$ is absolutely continuous with respect to $\bbQ^+$, so
\begin{equation} \label{EqFinalArgumentContradiction}
\bbQ^+\big(u_1(2)=u_2(2)\big) > 0.
\end{equation}
Since under $\bbQ^+$ the random variable $u_1(2)$ has law $\mu$ and the random variable $u_2(2)$ has law $\nu$ we cannot have at the same time \eqref{EqFinalArgumentContradiction} and the fact that $\mu$ and $\nu$ are singular. We thank M. Hairer for sharing his insight on this reasoning.   }
\end{Dem}

\ssk

Together with the existence result proved in \cite{BDFT}, Theorem \ref{ThmUniqueness} allows us to define the $\Phi^4_3$ measure on $M$ as the unique invariant probability measure of the semigroup $(\mcP_t)_{t\geq 0}$. This shows that the $\Phi^4_3$ measure is associated with the Riemannian manifold $M$. See the second paragraph after Theorem \ref{ThmMainInvariance} for the fact that one should more properly speak of a one parameter family of $\Phi^4_3$ measures.

\ssk

We denote below by $\mcL(X)$ the law of a random variable $X$.

\ssk

\begin{cor}
The total variation distance between the law of $u(t)$ and the invariant probability measure decays exponentially fast.
\end{cor}

\ssk

\begin{Dem}
The estimate \eqref{EqUniformBoundNonCoupling} is actually uniform with respect to the choice of initial conditions for $u_1$ and $u_2$. It follows classically from that fact and the Markovian character of the dynamics that the total variation distance between the law of $u(t)$ and the invariant probability measure decays exponentially fast
$$
\big\Vert \mcL(u_1(t)) - \mcL(u_2(t))\big\Vert_{\textsc{TV}} \leq (1-a)\big\Vert \mcL(u_1(t-1)) - \mcL(u_2(t-1))\big\Vert_{\textsc{TV}} \leq \cdots \lesssim (1-a)^t.
$$
\end{Dem}

\medskip

\section{Strong Feller property}
\label{SubsectionStrongFeller}

The coupling used in the proof Theorem \ref{ThmUniqueness} can be used to prove the strong Feller property of the semigroup $(\mcP_t)_{t\geq 0}$ by showing that it satisfies some Harnack-type inequality. {\color{black} As a preliminary remark to the next statement note that since one has the inclusion
$$
\big\{\tau(\ell,\phi_1,\phi_2)=2\big\} \subset \bigg\{\big\Vert v(1,\phi'_1) - v_\ell(1,\phi'_2) \big\Vert_{L^2} > m(\widehat\xi\,) \, \frac{\ell - o(\widehat\xi ; \ell)}{2} \bigg\}
$$
with $v_\ell(1,\phi'_2) = v(1,\phi'_2)$ it follows from the $L^p$ or $L^\infty$ coming down from infinity result that
$$
\bbP\big(\tau(\ell,\phi_1,\phi_2)=2\big) \eqdef {\color{black} a_\ell}= o_\ell(1)
$$
is a function of $\ell$ that goes to $0$ as $\ell$ goes to infinity.

\ssk

\begin{thm} \label{ThmHarnack}
Pick a finite exponent $p_1>1$ and a time $t>0$. For any $\ell>0$ there exists a function 
$$
\Psi_\ell : C^{-1/2-\epsilon}(M)\times C^{-1/2-\epsilon}(M) \rightarrow \bbR,
$$ 
that is null on the diagonal and continuous, such that the inequality
\begin{equation} \label{EqHarnack}
(\mcP_tf)^{p_1}(\phi_2) \leq \mcP_t(f^{p_1})(\phi_1) \, e^{\Psi_\ell(\phi_1,\phi_2)} \Big(1 + {\color{black} a_\ell}^{\frac{1}{p_1}}\Vert f\Vert_\infty \, e^{\Psi_\ell(\phi_1,\phi_2)} \Big)^{p_1},
\end{equation}
holds for any measurable bounded function $f\geq 1$ on $C^{-1/2-\epsilon}(M)$ and any $\phi_1,\phi_2\in C^{-1/2-\epsilon}(M)$.
\end{thm}

\ssk

\begin{Dem}
We use the notations of the proof of Theorem \ref{ThmUniqueness}. Since $u_{\ell}(2,\phi_2) = u(2,\phi_1)$ on the event $\{\tau(\ell,\phi_1,\phi_2) < 2\}$ one has
\begin{equation*} \begin{split}
(\mcP_t f)(\phi_2) &= \bbE_{\bbQ_{\ell,\phi_1,\phi_2}}\big[f(u_{\ell}(2,\phi_2))\big] = \bbE\big[R_{\ell,\phi_1,\phi_2}f(u_{\ell}(2,\phi_2))\big]   \\
&=\bbE\big[R_{\ell,\phi_1,\phi_2} f\big(u(2,\phi_1)\big) {\bf 1}_{\tau<2}\big] + \bbE_{\bbQ_{\ell,\phi_1,\phi_2}}\big[f\big(u_{\ell}(2,\phi_2)\big) {\bf 1}_{\tau=2}\big]
\end{split} \end{equation*}
and we obtain an inequality of the form \eqref{EqHarnack} from H\"older inequality and the condition $f\geq 1$, that allows to factorize in the second inequality below,

\begin{equation*} \begin{split}
(\mcP_t f)^{p_1}(\phi_2) &\leq \Big( \bbE\big[R_{\ell,\phi_1,\phi_2}  f(u(2,\phi_1))\big] + \Vert f\Vert_\infty \,\bbE\big[ R_{\ell(\alpha),\phi_1,\phi_2}^{\frac{p_1}{p_1-1}}\big]^{\frac{p_1-1}{p_1}} {\color{black} a_\ell}^{\frac{1}{p_1}} \Big)^{p_1}   \\
&\leq \bbE\big[R_{\ell,\phi_1,\phi_2} f(u(2,\phi_1))\big]^{p_1} \Big( 1 + \Vert f\Vert_\infty \,\bbE\big[ R_{\ell(\alpha),\phi_1,\phi_2}^{\frac{p_1}{p_1-1}}\big]^{\frac{p_1-1}{p_1}} {\color{black} a_\ell}^{\frac{1}{p_1}}  \Big)^{p_1}   \\
\end{split} \end{equation*}
\begin{equation*} \begin{split}
				    &\leq \bbE\big[f^{p_1}(u(2,\phi_1))\big] \, \bbE\big[ R_{\ell,\phi_1,\phi_2}^{\frac{p_1}{p_1-1}}\big]^{p_1-1} \Big( 1 + \Vert f\Vert_\infty \,\bbE\big[ R_{\ell(\alpha),\phi_1,\phi_2}^{\frac{p_1}{p_1-1}}\big]^{\frac{p_1-1}{p_1}} {\color{black} a_\ell}^{\frac{1}{p_1}} \Big)^{p_1}.
\end{split} \end{equation*}
This is \eqref{EqHarnack} with 
$$
e^{\Psi_\ell(\phi_1,\phi_2)} = \bbE\big[ R_{\ell,\phi_1,\phi_2}^{\frac{p_1}{p_1-1}}\big]^{p_1-1}
$$
(The function $\Psi_\ell$ also depends on $p_1$ but we do not emphasize that dependence in the notation as it is irrelevant for us here.) We check from classical arguments that $\Psi_\ell(\phi_1,\phi_2)$ is a continuous function of $\phi_1$ and $\phi_2$.
\end{Dem}

\ssk

\begin{cor}
The semigroup $(\mcP_t)_{t\geq 0}$ has the strong Feller property. 
\end{cor}

\ssk

\begin{Dem}
We follow F.Y. Wang's classical proof -- see e.g. Theorem 1.4.1 in \cite{WangBook}. Fix $t>0$. Applying \eqref{EqHarnack} to $f=1+r g$, with a measurable function $0\leq g\leq 1$ and $0<r\leq 1$, one gets
\begin{equation} \label{EqTool}
\big(1+r (\mcP_t g)(\phi_2)\big)^{p_1} \leq \mcP_t\big((1+r g)^{p_1}\big)(\phi_1) \, e^{\Psi_\ell(\phi_1,\phi_2)} \Big( 1 + {\color{black} a_\ell}^{\frac{1}{p_1}} (1+\Vert g\Vert_\infty) \, e^{\Psi_\ell(\phi_1,\phi_2)} \Big)^{p_1}
\end{equation}
so
\begin{equation*} \begin{split}
1 &+ p_1r (\mcP_t g)(\phi_2) + o(r)   \\
&\leq \Big(1+p_1r(\mcP_tg)(\phi_1) + o(r)\Big) e^{\Psi_\ell(\phi_1,\phi_2)} \Big( 1 + {\color{black} a_\ell}^{\frac{1}{p_1}}(1+\Vert g\Vert_\infty) \, e^{\Psi_\ell(\phi_1,\phi_2)} \Big)^{p_1}.
\end{split} \end{equation*}
Sending $\phi_2$ to $\phi_1$ one gets from the continuity of $\Psi_\ell$
$$
1 + p_1r \, \underset{\phi_2\rightarrow\phi_1}{\limsup} \, (\mcP_t g)(\phi_2) + o(r) \leq \Big(1+p_1r (\mcP_tg)(\phi_1) + o(r)\Big) \Big( 1 + {\color{black} a_\ell}^{\frac{1}{p_1}} (1+\Vert g\Vert_\infty)\Big)^{p_1}
$$
and since $\ell>0$ and $r>0$ are arbitrary {\color{black} and ${\color{black} a_\ell}$ goes to $0$ as $\ell$ goes to infinity}
$$
\underset{\phi_2\rightarrow\phi_1}{\limsup} \, (\mcP_t g)(\phi_2) \leq (\mcP_t g)(\phi_1).
$$
Exchanging $\phi_1$ and $\phi_2$ in \eqref{EqTool}, the same reasoning gives
$$
(\mcP_t g)(\phi_1) \leq \underset{\phi_2\rightarrow\phi_1}{\liminf} \, (\mcP_t g)(\phi_2),
$$
from which the continuity of $\mcP_t g$ at $\phi_1$ follows. The conclusion follows since $\phi_1$ is arbitrary.   \vspace{0.15cm}
\end{Dem}

\medskip

The remainder of this work is dedicated to proving the a priori estimates that were involved in the assumptions of Lemma \ref{LemMonotony}.

\medskip

\section{An $L^\infty$ coming down from infinity result}
\label{SectionLinftyComingDown}

We proved in Section 2 of \cite{BDFT} an explicit control on the $L^p$ norm ($1\leq p<\infty$) of the solution $v$ to Equation \eqref{EqJPFormulationSynthetic} that is independent of its initial condition. However this $L^p$ coming down from infinity result is not sufficient for our needs here. We obtain in the present section an $L^\infty$ coming down from infinity result for $v$ that is the content of Theorem \ref{ThmLinftyComingDown} below. Its proof will follow closely the seminal work \cite{MoinatWeber} of Moinat \& Weber on the $\Phi^4_3$ equation \eqref{EqPhi43SPDE} itself.

\ssk

We need a piece of notation to state Theorem \ref{ThmLinftyComingDown} below. Throughout this section we will use the shorthand notation $\Vert\cdot\Vert$ for $\Vert\cdot\Vert_{L^\infty}$ and $\Vert\cdot\Vert_D$ for $\Vert\cdot\Vert_{L^\infty(D)}$, for a parabolic domain $D$, meaning a subset of $(0,T)\times M$. Set
$$
\mcT \defeq \Big\{ A, B, Z_2, Z_1, Z_0 \Big\}
$$
and 
\begin{equation} \label{EqChoiceExponents} \begin{split}
&\frac{1}{m_A} = \epsilon\Big(\frac{1}{2}-\epsilon\Big),\quad \frac{1}{m_B} = 1-\epsilon\Big(\frac{1}{2}-2\epsilon\Big) - 3\epsilon,   \\
&\frac{1}{m_{Z_2}} = \frac{1}{2} -\epsilon', \quad \frac{1}{m_{Z_1}} = \frac{3}{2}-\epsilon', \quad \frac{1}{m_{Z_0}} = \frac{5}{2}-\epsilon', 
\end{split} \end{equation}
with $\epsilon'>3\epsilon$. Fix $T\geq 2$; its precise value does not matter here. For $0<s<\sqrt{T}$ we define the parabolic domain
$$
\mcD_s \defeq  (s^2, T) \times M \subset \bbR\times M.
$$
For $\tau\in\mcT$ we define $[\tau]_{\vert\tau\vert}$ as the norm of $\tau\in C_TC^{\vert\tau\vert}(M)$, where 
$$
\vert A\vert = 1-\epsilon,  \quad  \vert B\vert=-\epsilon,  \quad  \vert Z_i\vert=-1/2-\epsilon.
$$ 
Set 
$$
A_+ \defeq \sup_{\mcD_0} A,  \quad  A_- \defeq \min_{\mcD_0} A
$$
and
$$
c_A \defeq \big(1+ \max( A_+ , A_-^{-1})\big)^2.
$$

\ssk

\begin{thm} \label{ThmLinftyComingDown}
There exists a positive constant $C$ such that any solution $v$ of Equation \eqref{EqJPFormulationSynthetic} satisfies for all $0<s\leq 1$ the estimate
\begin{equation} \label{EqMWBound}
\Vert v\Vert_{\mcD_s} \leq C \max\bigg\{\frac{1+A_-^{-1/2}}{s} , \big((c_A[\tau]_{\vert\tau\vert})^{m_\tau}\big)_{\tau\in\mcT}\bigg\}.
\end{equation}
\end{thm}

\ssk

\subsection{Tools for the proof$\boldmath{.}$ \hspace{0.1cm}}

We collect in this section two ingredients that will play a key role in the proof of Theorem \ref{ThmLinftyComingDown}: A Schauder type estimate and a corollary of the maximum principle. First we set the scene. We use some partition of unity ${\bf 1}_M = \sum_{i=1} ^{n_0} \chi_i$ of $M$ associated with some charts $\chi_i : U_i\rightarrow \bbR^3$ to pull back each piece $(\textrm{id}\times \chi_i)_\sharp\Lambda$ of a distribution $\Lambda$ on $\bbR\times M$ into $\bbR\times\bbR^3$ and define for $\alpha<0$ the local $\alpha$-Besov-H\"older norm on an open set $O\subset\bbR\times M$ as
$$
\Vert \Lambda\Vert_\alpha \defeq \sum_{i=1}^{n_0} \sup_{0<\delta\leq 1}\delta^{-\alpha} \big\Vert \big((\textrm{id}\times\chi_i)_\sharp\Lambda\big)_\delta\big\Vert_{L^{\infty}(\chi_i^{-1}(O))},
$$
where $(\cdot)_\delta$ is the parabolic regularization operator from Section 2 of Moinat \& Weber's work \cite{MoinatWeber}. For $0<\alpha<1$ and a domain $D\subset \bbR\times M$ we define the H\"older seminorm
$$
[h]_{\alpha,D} \defeq \underset{z\neq z' \in D}{\sup} \, \frac{\vert h(z)-h(z')\vert}{\vert z-z'\vert^\alpha}
$$
and for $1<\alpha<2$ we set
$$
[h]_{\alpha,D} \defeq \sup_{z'=(t',x')\neq z=(t,x)\in D} \sup_{\theta\in T_xM} \frac{\big| h(z')-h(z) - d(x',x) \, dh(z)(\theta)\big|}{d(z',z)^\alpha}.
$$
We write $[h]_{\alpha}$ for $[h]_{\alpha,\bbR\times M}$. For a function $h\in\mcC^\beta(\bbR\times M)$, with $0<\beta<1$, and any open set $O\subset \bbR\times M$ with $\delta$-parabolic neighbourhood $O_\delta$ we have
$$
\Vert h_\delta-h\Vert_{L^\infty(O)} \lesssim \delta^\beta [h]_{\beta,O_\delta},
$$
for an implicit multiplicative constant independent of $h$. For $z\in\bbR\times M$ we denote by $B(z,\delta)\subset \bbR\times M$ the parabolic ball of center $z$ and radius $\delta$. The next statement is a consequence of the reconstruction theorem in the case of a Young product. See e.g. Theorem 14.1 of Caravenna \& Zambotti's review \cite{CZ} and Theorem 18 of \cite{RinaldiSclavi} for its extension to a manifold setting. See also Moinat \& Weber's version of the reconstruction theorem, Theorem 2.8 in \cite{MoinatWeber}.

\ssk

\begin{lem}\label{LemCommutatorEstimate} 
Pick $\alpha <0$ and $\beta\in (0,1)$ such that $\alpha+\beta>0$.  Let $f\in C^\alpha(B(z, 2\delta))$ and $g\in C^\beta(B(z, \delta))$. Then we have
\begin{equation}
\big|\big( (fg)_\delta -f_\delta g\big)(z) \big| \lesssim  \delta^{\alpha+\beta} [f]_{\alpha, B(z, \delta)} \, [g]_{\beta, B(z, \delta)}.
\end{equation}
Moreover if $f\in L^\infty  (B(z, \delta)) $ we have
\begin{equation} \label{EqCommutatorEstimateRegularization}
\big| (fg)_\delta(z) - (f_\delta g)(z) \big| \lesssim  \delta^{\beta} \|f\|_{\infty} \, [g]_{\beta, B(z, \delta)}.
\end{equation}
\end{lem}

\ssk

The following notation is involved in the inequalities \eqref{eqDVUniformSchauder} and \eqref{eqDVIncrementSchauder} below. For $0<\gamma<1$ and for any continuous function $W$ on $\mcD_s^2$ and subset $\mcE$ of $\mcD_s$ we set for $\rho>0$
$$
\Vert W\Vert_{(\gamma\,;\,\rho), \mcE} \defeq \underset{z'\neq z\in \mcE, \vert z'-z\vert\leq \rho}{\sup}\,\frac{\vert W(z',z)\vert}{\vert z'-z\vert^\gamma};
$$
this is a kind of $\rho$-local $\gamma$-H\"older norm of $W$ on the set $\mcE$. For a function $w$ on $\mcD_s$ we denote by $\Vert w\Vert_{(\gamma\,;\,\rho), \mcE}$ the $\Vert \cdot\Vert_{(\gamma\,;\,\rho), \mcE}$-norm of the two variable function $w(z')-w(z)$. The quantity $\Vert w\Vert_{(0\,;\,\rho), \mcE}$ is for instance the oscillation of $w$ in balls of $\mcE$ of radius $\rho$. The following statement provides a strong control on a function $w$ in terms of a control on $(\partial_t-\Delta) w$ and a `weak' control of $w$ itself. It is a simplified version of a subtler and finer estimate proved by Moinat \& Weber in Lemma 2.11 of \cite{MoinatWeber}. The latter was itself a generalization of Proposition 2 of Otto, Sauer, Smith \& Weber's work \cite{OSSW}.

\ssk

\begin{thm} \label{LemLocalSchauder} 
Let a regularity exponent $\kappa\in (1,2)$ and a constant $\delta_0>0$ be given. There is a constant $c_1>0$ with the following property. If for all $0<4\delta\leq \lambda \leq \lambda_0$ one has
\begin{eqnarray} \label{EqAssumptionSchauder}
\delta^{2-\kappa} \big\Vert \big((\partial_t-\Delta+1) w\big)_\delta \big\Vert_{\mcD_{s+\lambda-\delta}}\leq C_\lambda 
\end{eqnarray}
for some function $w : \mcD_s\rightarrow\bbR$, then one has
\begin{equation} \label{EqSchauderMain}
\sup_{0<\lambda\leq \lambda_0}\lambda^\kappa [w]_{\kappa,\mcD_{s+\lambda}} \leq {\color{black} c_1} \Big(\sup_{\lambda\leq \lambda_0} \lambda^{\kappa} C_\lambda + \Vert w \Vert_{\mcD_s}\Big).
\end{equation}
Moreover one can associate to any $0<\delta<\delta_0$ a constant $\rho_\delta>0$ such that for $0<\rho<\rho_\delta$ one has the following form of Schauder estimate
\begin{equation} \label{eqDVUniformSchauder}
\Vert dw\Vert_{\mcD_{s+\delta}} \lesssim \rho^{\kappa-1} [w]_{\kappa,\mcD_{s+\delta}} + \frac{1}{\rho}\,\Vert w\Vert_{(0\,;\,\rho), \mcD_{s+\delta}} 
\end{equation}
and 
\begin{equation} \label{eqDVIncrementSchauder}
 [dw]_{\kappa-1,\mcD_{s+\delta}} \lesssim [w]_{\kappa,\mcD_{s+\delta}} + \frac{1}{\rho^\kappa}\, \Vert w\Vert_{(0\,;\,\rho),\mcD_{s+\delta}}.
\end{equation}
\end{thm}

\ssk

In addition to Theorem \ref{LemLocalSchauder} we will also use in the proof of Theorem \ref{ThmLinftyComingDown} the following statement which explains the form of the bound \eqref{EqMWBound}.

\medskip

\begin{lem} \label{LemMaximumPrinciple}
Let $f,g,h$ be some continuous functions on $[0,1]\times M$ with $\min g\eqdef g_- > 0$. Any continuous function $w$ on $[0,1]\times M$ such that
\begin{equation} \label{EqContinuousFormMaxPrinciple}
(\partial_t-\Delta + f\cdot\nabla)w = - g w^3 + h
\end{equation}
on $(0,1)\times M$ satisfies for all $0<t\leq 1$ the estimate
$$
\Vert w\Vert_{[t,1]\times M} \leq  \max\big((g_-\,t)^{-\frac{1}{2}}, (\|h\|_\infty/g_-)^{\frac{1}{3}}\big).
$$
\end{lem}

\medskip

Indeed one can check that the functions
$$
\pm\Big((2g_-t)^{-\frac{1}{2}} + (\|h\|_\infty/g_- )^{\frac{1}{3}}\Big)
$$
are supersolution (+) and subsolution (-) of Equation \eqref{EqContinuousFormMaxPrinciple}, so the conclusion comes from the comparison/maximum principle. (See for instance Section 5.2 of Taylor's textbook \cite{TaylorPDE}.) The strong damping effect of the superlinear term $-Av^3$ in Equation \eqref{EqJPFormulationSynthetic} will only be used in the proof of Theorem \ref{ThmLinftyComingDown} by appealing to Lemma \ref{LemMaximumPrinciple} on a regularized version of Equation \eqref{EqJPFormulationSynthetic}.

\medskip

\subsection{Proof of Theorem \ref{ThmLinftyComingDown}$\boldmath{.}$ \hspace{0.1cm}}
\label{SectionProofLInftyComingDown}

The main step of the proof of Theorem \ref{ThmLinftyComingDown} consists in showing that if one has
$$
{\color{black} \|v\|_{\mcD_s} \geq 32}
$$
and
$$
{\color{black} [Z_1-1]_{-1/2-\epsilon} \leq \frac{c}{c_A} \, \|v\|^{1/m_{Z_1}}_{\mcD_s} }
$$
and the inequalities
\begin{equation} \label{EqAssumptionMWProof}
[\tau]_{\vert\tau\vert} \leq \frac{c}{c_A} \, \|v\|^{1/m_\tau}_{\mcD_s} \qquad \big(\tau\in \big\{A,B, Z_2, Z_0\big\}\big)
\end{equation}
for a well-chosen fixed positive constant $c$ independent of $v$ then
\begin{equation} \label{EqMainStepProofComingDown}
\Vert v\Vert_{\mcD_{s+s_1}} \leq \max\bigg\{\frac{2 A_-^{-\frac{1}{2}}}{s_1}, \frac{\Vert  v\Vert_{\mcD_s}}{2}\bigg\}
\end{equation}
for all $s_1$ with $s+s_1\leq 1/2$ and $s_1\geq 1/\|v\|_{\mcD_s}$. The proof of this inequality is the content of item {\it (a)} below. We explain in item {\it (b)} how the statement of Theorem \ref{ThmLinftyComingDown} follows from that fact. 

\medskip

\noindent {\it (a) The main step: Proof of \eqref{EqMainStepProofComingDown}.} The proof of \eqref{EqMainStepProofComingDown} proceeds in three steps. 
\begin{itemize}
	\item[\textbf{\textsf{1.}}] We prove that $\Vert  v\Vert_{\mcD_s}$ controls the $(3/2-\epsilon)$ and $(1/2\pm\epsilon)$ seminorms of $v$ on the smaller parabolic domains $\mcD_{s+\lambda}$. (The Schauder estimate from Theorem \ref{LemLocalSchauder} is used for that purpose.)
	\item[\textbf{\textsf{2.}}] We apply Lemma \ref{LemMaximumPrinciple} to a $\delta$-regularized version of Equation \eqref{EqJPFormulationSynthetic}. Together with the result of Step 1 it provides a uniform bound on $v$ on domains of the form $\mcD_{s+s_1}$. Both $s_1$ and the regularization parameter are free in that step.
	\item[\textbf{\textsf{3.}}] We tune the regularization parameter $\delta$ to optimize the bound from Step 2.
\end{itemize}
We use below the shorthand notation $\mcL$ for the differential operator $\partial_t-\Delta+1$. The constant $c$ in \eqref{EqAssumptionMWProof} will be chosen later, just before \eqref{eq:l_infty_main_est}, and we write 
$$
c'_A \defeq \frac{c}{c_A}.
$$

\ssk

\noindent $\rhd$ \textbf{\textsf{Step 1.}} {\it We first derive from Equation \eqref{EqJPFormulationSynthetic} and the Schauder estimate from Theorem \ref{LemLocalSchauder} some $\lambda$-dependent bound on the $(3/2-\epsilon)$-H\"older norm of $v$ on $\mcD_{s+\lambda}$ in terms of its uniform norm on the larger domain $\mcD_s$.} (The bound explodes as $\lambda$ goes to $0$.) We start from the regularized version 
\begin{equation} \label{LVDelta0}
(\mcL v)_\delta = - (Av^3)_\delta + (B\cdot \nabla v)_\delta + (Z_2 v^2)_\delta + (Z_1v)_\delta + (Z_0)_\delta. 
\end{equation}
of Equation \ref{EqJPFormulationSynthetic}. Note that the conclusion of Theorem \ref{LemLocalSchauder} makes it possible to use the $(3/2-\epsilon)$-H\"older seminorm of $v$ in our estimates of the right-hand side of \eqref{LVDelta0} provided it comes with a small factor that can eventually be absorbed in the left-hand side of \eqref{EqSchauderMain}. For the $Z_0$ and $A$ terms in \eqref{LVDelta0} we simply bound
\begin{equation*} \begin{split}
\| (Z_0)_\delta\|_{\mcD_{s+\lambda-\delta}} &\leq c'_A \, \delta^{-\frac{1}{2} - \epsilon} \|v\|_{\mcD_s}^{1/m_{Z_0}}   \\
\| (Av^3)_\delta \|_{\mcD_{s+\lambda-\delta}} &\leq \sqrt{c_A} \, \|v\|_{\mcD_s}^{3},  \quad \text{ since }  \|A\|_{\mcD_0}\leq  \sqrt{c_A}.
\end{split} \end{equation*}
For the $Z_2$ and $Z_1$ terms in \eqref{LVDelta0} one gets from the assumption \eqref{EqAssumptionMWProof} for $\delta\leq \lambda/4$ 

\begin{equation*} \begin{split}
\|(Z_2v^2)_{\delta}\|_{\mcD_{s+\lambda-\delta}} &\leq c'_A\delta^{-\frac{1}{2}-\epsilon}\|v\|_{\mcD_s}^{2+1/m_{Z_2}} + 2 c'_A \delta^{-\epsilon} [v]_{(1/2+2\epsilon\,;\,\delta), \mcD_{s +\lambda/2}}   \|v\|_{\mcD_s}^{1/m_{Z_2}+ 1},   \\
\|(Z_1v)_{\delta}\|_{\mcD_{s+\lambda- \delta}} &\leq c'_A\delta^{-\frac{1}{2}-\epsilon |}\|v\|_{\mcD_s}^{1+1/m_{Z_1}} + c'_A \delta^{-\epsilon} [v]_{(1/2+2\epsilon\,;\,\delta), \mcD_{s+\lambda/2}}  \|v\|_{\mcD_s}^{1/m_{Z_1}} .
\end{split} \end{equation*}
We used here an upper bound for the $C_TC^{1/2+2\epsilon}(M)$ norm of $v$ to consider the product $Z_2v$ as a well-defined Young product, a continuous function of $Z_2$ and $v$. It turns out to be useful to introduce the commutator $B_\delta\cdot\nabla v - (B\cdot\nabla v)_\delta$ to estimate $(B\cdot\nabla v)_\delta$ itself as we get from Lemma \ref{LemCommutatorEstimate} the bound
\begin{equation*} \begin{split}
\Vert (B\cdot\nabla v)_\delta\Vert_{\mcD_{s+\lambda-\delta}} &\leq \big\Vert B_\delta\cdot\nabla v - (B\cdot\nabla v)_\delta \big\Vert_{\mcD_{s+\lambda-\delta}} + \Vert B_\delta\cdot\nabla v\Vert_{\mcD_{s+\lambda-\delta}}   \\
&\lesssim \delta^{1/2-2\epsilon} [B]_{-\epsilon} [\nabla v]_{\frac{1}{2}-\epsilon, \mcD_{s+\lambda/2}} + \delta^{-\epsilon} [B]_{-\epsilon} \Vert \nabla v\Vert_{\mcD_{s+\lambda/2}}.
\end{split} \end{equation*}
Note that for all $\lambda$ sufficiently small one has from \eqref{eqDVUniformSchauder} and \eqref{eqDVIncrementSchauder} 
\begin{equation} \label{EqGradientBoundUniform}
\Vert \nabla v\Vert_{\mcD_{s+\lambda/2}} \lesssim \lambda^{\frac{1}{2}-\epsilon} [v]_{3/2-\epsilon, \mcD_{s+\lambda/2}} + \lambda^{-1} \Vert v\Vert_{(0;\lambda), \mcD_{s+\lambda/2}},
\end{equation}
and 
\begin{equation} \label{EqGradientBound}
[\nabla v]_{1/2-\epsilon, \mcD_{s+\lambda/2}} \lesssim [v]_{3/2-\epsilon, \mcD_{s+\lambda/2}} + \lambda^{-\frac{3}{2}+\epsilon} \Vert v\Vert_{(0;\lambda), \mcD_{s+\lambda/2}}.
\end{equation}
Combining these estimates all together yields the bound
$$
\delta^{1/2+\epsilon} \|(\mcL v)_\delta \|_{\mcD _{s+\lambda-\delta}} \lesssim C_\lambda
$$
where 
\begin{equation*} \begin{split}
C_\lambda &\defeq \lambda^{\frac{1}{2}+\epsilon} \sqrt{c_A} \, \|v\|_{\mcD_s}^{3} + c'_A \lambda \|v\|_{\mcD_s}^{1/m_B} [v]_{3/2-\epsilon, \mcD_{s+\lambda/2}} + c'_A\lambda^{-\frac{1}{2}} \|v\|^{1+1/m_B}_{\mcD_s}   \\
&\quad+ c'_A \|v\|_{\mcD_s}^{2+1/m_{Z_2}} + c'_A \lambda^{\frac{1}{2}} [v]_{(1/2+2\epsilon\,;\,\delta), \mcD_{s+\lambda/2}}   \|v\|_{\mcD_s}^{1/m_{Z_2}+ 1}   \\
&\quad+ c'_A \|v\|_{\mcD_s}^{1+1/m_{Z_1}} + c'_A \lambda^{\frac{1}{2}} [v]_{(1/2+2\epsilon\,;\,\delta), \mcD_{s+\lambda/2}}   \|v\|_{\mcD_s}^{1/m_{Z_1}} + c'_A  \|v\|_{\mcD_s}^{1/m_{Z_0}}.  
\end{split} \end{equation*}
Note that $\mcD_{s+\lambda/2}$ is involved in the definition of $C_\lambda$ rather than $\mcD_{s+\lambda}$. The Schauder estimate from Theorem \ref{LemLocalSchauder} gives us the second inequality below
\begin{equation} \label{EqAlmostThereStep1}
\underset{\lambda\leq\frac{\lambda_0}{2}}{\sup}\; \lambda^{\frac{3}{2}-\epsilon} [u]_{\frac{3}{2}-\epsilon,\mcD_{s+\lambda}} \leq \underset{\lambda\leq\lambda_0}{\sup}\; \lambda^{\frac{3}{2}-\epsilon} [u]_{3/2-\epsilon,\mcD_{s+\lambda}} \lesssim \underset{\lambda\leq\lambda_0}{\sup}\; \lambda^{\frac{3}{2}-\epsilon} C_\lambda + \Vert v\Vert_{\mcD_s},
\end{equation}
where the same domains are involved in the supremum on both sides. Even though $\lambda^{3/2-\epsilon} C_\lambda$ depends on the $(3/2-\epsilon)$-seminorm of $v$ on $\mcD_{s+\lambda/2}$ it comes with a factor that will eventually be small for the choice of $\lambda_0$ made below in \eqref{EqChoiceLambda0}. Recall $0<\delta\leq\lambda/4$. The constant $\lambda^{3/2-\epsilon} C_\lambda$ will still depend on the $\delta$-local $(1/2-\epsilon)$-seminorm of $v$ on $\mcD_{s+\lambda/2}$. To eventually have a bound on the supremum in the right hand side of \eqref{EqAlmostThereStep1} that only involves $\|v\|_{\mcD_s}$ we use the elementary estimate
\begin{equation} \label{EqNorm1sur2v} 
[v]_{(1/2+2\epsilon;\delta),\mcD_{s+\lambda}} \leq \lambda^{1-3\epsilon} [v]_{3/2-\epsilon,\mcD_{s+\lambda}} + \lambda^{\frac{1}{2} - 2\epsilon} \Vert \nabla v\Vert_{\mcD_{s+\lambda}}
\end{equation}
together with \eqref{EqGradientBoundUniform} to see that 
$$
[v]_{(1/2+2\epsilon;\delta),\mcD_{s+\lambda}} \lesssim \lambda^{1-3\epsilon} [v]_{3/2-\epsilon, \mcD_{s+\lambda}} + \lambda^{-\frac{1}{2} - 2\epsilon} \Vert v\Vert_{\mcD_{s+\lambda}}.
$$
Recall $\epsilon' > 3\epsilon$, where $\epsilon'$ is involved in the definition of the exponents $m_{Z_i}$ in \eqref{EqChoiceExponents}. Overall we have a $[v]_{3/2-\epsilon,\mcD_{s+\lambda}}$ term in an upper bound for $\lambda^{3/2-\epsilon} C_\lambda$ that comes with a factor 
\begin{equation} \label{EqFactor}
\lambda \|v\|_{\mcD_s}^{1/m_B} + \lambda^{3/2-3\epsilon} \big(\|v\|_{\mcD_s}^{1+1/m_{Z_2}} + \|v\|_{\mcD_s}^{1/m_{Z_1}}\big) = \|v\|_{\mcD_s}^{-\frac{7\epsilon}{2}+2\epsilon^2} + 2\lambda^{3/2-3\epsilon} \|v\|_{\mcD_s}^{3/2+\epsilon'}.
\end{equation}
Choosing
\begin{equation} \label{EqChoiceLambda0}
\lambda_0 \leq \Vert v\Vert_{\mcD_s}^{-1}
\end{equation}
ensures that the factor \eqref{EqFactor} is small, so the corresponding term can be absorbed in the left hand side of \eqref{EqAlmostThereStep1} and we get
\begin{equation} \label{eq:est_3_2}
\sup_{\lambda\leq \lambda_0/2} \lambda^{\frac{3}{2}- \epsilon} [v]_{3/2-\epsilon, \mcD_{s+\lambda}}\lesssim \|v\|_{\mcD_s}.  
\end{equation}
Then it follows from \eqref{EqNorm1sur2v} that 
\begin{equation} \label{eq:est_1_2}
\sup_{\lambda\leq \lambda_0/2} \lambda^{\frac{1}{2}+2\epsilon} [v]_{(\frac{1}{2}+2\epsilon; \lambda/4), {\mcD_{s+\lambda}}} \lesssim \|v\|_{\mcD_s}.
\end{equation}
A similar estimate holds for $\sup_{\lambda\leq \lambda_0/2} \lambda^{\frac{1}{2}-\epsilon} [v]_{(\frac{1}{2}-\epsilon; \lambda/4), {\mcD_{s+\lambda}}}$.

\medskip

\noindent $\rhd$ \textbf{\textsf{Step 2.}} \textit{We regularize Equation \eqref{EqJPFormulationSynthetic} and apply the maximum principle of Lemma \ref{LemMaximumPrinciple} to its solution to obtain an upper bound on $\Vert v_\delta\Vert_{\mcD_{s+s_1}}$, for $s+s_1<1/2$ and $0<\delta\ll s_1$, in terms of $\Vert v\Vert_{\mcD_s}$ and $\delta$. (The inequality $\delta\ll s_1$ is quantified just before \eqref{EqBoundsStep2AB} below).} The regularized version of Equation \eqref{EqJPFormulationSynthetic} takes the form
\begin{equation*} \begin{split}
(\partial_t-\Delta + B_\delta\cdot \nabla) v_\delta &= - A v_\delta^3 + [\mcL , (\cdot)_\delta](v) + \big( B_\delta\cdot \nabla v_\delta - (B\cdot \nabla v)_\delta \big)   \\
&\quad+ \big(Av_\delta^3 - (Av^3)_\delta\big) + (Z_2 v^2)_\delta + \big((Z_1-1)v\big)_\delta + (Z_0)_\delta.  \\
&\eqdef -Av_\delta^3 + h_\delta.
\end{split} \end{equation*}
(The $-1$ in the linear term in $v$ come from the fact that $\mcL=\partial_t-\Delta+1$ while Lemma \ref{LemMaximumPrinciple} involves the operator $\partial_t-\Delta$.) For all $s_1>0$ such that $s+s_1<\frac{1}{2}$, the pointwise estimate from Lemma \ref{LemMaximumPrinciple} gives here for $\|v_\delta\|_{\mcD_{s+s_1}} $ the upper bound
\begin{equation} \label{EqEstimateStep2} \begin{split}
\max&\Big\{ \big(A_-^{-\frac{1}{2}} \, \frac{2}{s_1} \,,\, A_-^{-\frac{1}{3}} \big\| [\mcL , (\cdot)_\delta](v) \big\|_{\mcD_{s+s_1/2}}^{\frac{1}{3}} \, , \, A_-^{-\frac{1}{3}} \| A v_\delta^3- (Av^3)_\delta\|_{\mcD_{s+s_1/2}}^{\frac{1}{3}}   \\
&\quad A_-^{-\frac{1}{3}}\|B_\delta\cdot \nabla v_\delta -(B\cdot \nabla v)_\delta \|_{\mcD_{s+s_1/2}}^{\frac{1}{3}} \,,\, A_-^{-\frac{1}{3}} \| (Z_2 v^2)_\delta \|_{\mcD_{s+s_1/2}}^{\frac{1}{3}}, A_-^{-\frac{1}{3}} \| (Z_1 v)_\delta \|_{\mcD_{s+s_1/2}}^{\frac{1}{3}},  \\
&\quad A_-^{-\frac{1}{3}} \| (Z_0 )_\delta \|_{\mcD_{s+s_1/2}}^{\frac{1}{3}}  \Big\},
\end{split} \end{equation}
and one gets a bound on $v$ writing for $0<\delta\leq s_1$
\begin{equation} \label{EqVFromVdelta}
\|v\|_{\mcD_{s+s_1}} \leq \|v_\delta\|_{\mcD_{s+s_1}} + \delta^{\frac{1}{2}-\epsilon} [v]_{(\frac{1}{2}-\epsilon\,;\,2\delta), \mcD_{s+s_1-\delta}} \lesssim   \|v_\delta\|_{\mcD_{s+s_1}} + \Big(\frac{\delta}{s_1}\Big)^{\frac{1}{2}-\epsilon} \Vert v\Vert_{\mcD_s},
\end{equation}
as a consequence of the reconstruction estimate \eqref{EqCommutatorEstimateRegularization} for the first inequality and \eqref{eq:est_1_2} for the second inequality. Let us introduce a positive parameter $k\geq 4$ that will be chosen in \eqref{EqChoiceK} below in Step 3. Leaving aside the exponent $1/3$, the different terms in the above maximum can be bounded as follows for $0<\delta\leq s_1/k$
\begin{equation} \label{EqBoundsStep2AB} \begin{split}
&\big\| [\mcL , (\cdot)_\delta](v) \big\|_{\mcD_{s+s_1/2}} \lesssim \delta^{-1} \Vert v\Vert_{\mcD_s},   \\
&\big\| A v_\delta^3 - (Av^3)_\delta \big\|_{\mcD_{s+s_1/2}} \lesssim \Big(\delta^{\frac{1}{2}}\Vert A\Vert_{1/2, \mcD_s} + k^{-\frac{1}{2}+\epsilon} \Vert A\Vert_{\mcD_s}\Big) \Vert v\Vert_{\mcD_s}^3,   \\
&\big\| B_\delta \cdot \nabla v_\delta - (B\cdot \nabla v )_\delta \big\|_{\mcD_{s+s_1/2}} \lesssim [B]_{-\epsilon} \Vert v\Vert_{\mcD_s} k^{-\frac{1}{2}+2\epsilon} s_1^{-1+3\epsilon},   \\
\end{split} \end{equation}
and
\begin{equation} \label{EqBoundsStep2Z} \begin{split}
&\| (Z_2 v^2)_\delta \|_{\mcD_{s+s_1/2}} \lesssim [Z_2]_{-1/2-\epsilon} \|v\|^2_{\mcD_s} \delta^{-\frac{1}{2}-\epsilon},   \\
&\| (Z_1 v)_\delta \|_{\mcD_{s+s_1/2}} \lesssim \delta^{-\frac{1}{2}-\epsilon} [Z_1]_{-1/2-\epsilon} \Vert v\Vert_{\mcD_s},   \\
&\| (Z_0)_\delta \|_{\mcD_{s+s_1/2}} \lesssim \delta^{-\frac{1}{2}-\epsilon} [Z_0]_{-1/2-\epsilon}.
\end{split} \end{equation}
$\bullet$ Indeed for the $A$ term one has
\begin{equation*} \begin{split}
(Av_\delta^3)(z) &- (Av^3)_\delta(z)   \\
&= \int \varphi_\delta(z,z')v(z')\Big\{\big(A(z)-A(z')\big)v_\delta^2(z) + A(z')\big(v^2_\delta(z)-v^2(z')\big)\Big\}dz'.
\end{split} \end{equation*}
Denote by $w_A(\cdot,D)$ the modulus of continuity of $A$ on a domain $D$. Above, the term with the increment of $A$ gives a contribution bounded above by $w_A(\delta , \mcD_{s+s_1/4}) \Vert v\Vert_{\mcD_s}^3$. As $0<\delta\leq s_1/k$ one has for $z\in\mcD_{s+s_1/2}$ and $z'\in\mcD_{s+s_1/4}$
$$
\big\vert v^2_\delta(z) - v^2(z')\big\vert \lesssim \delta^{1/2-\epsilon} [v]_{1/2-\epsilon, \mcD_{s+s_1/4}} \Vert v\Vert_{\mcD_s}, 
$$
and we get from the fact that $A$ is (better than) $1/2$-H\"older and the bound \eqref{eq:est_1_2} the estimate 
$$
\big\vert (Av_\delta^3)(z)  - (Av^3)_\delta(z) \big\vert \lesssim \Big(\delta^{\frac{1}{2}}\Vert A\Vert_{1/2, \mcD_s} + \Big(\frac{\delta}{s_1}\Big)^{\frac{1}{2}-\epsilon} \Vert A\Vert_{\mcD_s}\Big) \Vert v\Vert_{\mcD_s}^3.
$$
The condition $0<\delta\leq s_1/k$ gives the estimate for $A$ from \eqref{EqBoundsStep2AB}.

\noindent $\bullet$ For the $B$ term write
$$
B_\delta\cdot\nabla v_\delta - (B\cdot\nabla v)_\delta = \big(B_\delta\cdot\nabla v - (B\cdot\nabla v)_\delta\big) + B_\delta\cdot\nabla (v-v_\delta).
$$
We use Lemma \ref{LemCommutatorEstimate} to estimate the term in the big parenthesis on the right hand side. This gives
$$
\big\Vert B_\delta\cdot\nabla v_\delta - (B\cdot\nabla v)_\delta \big\Vert_{\mcD_s+s_1/2} \lesssim \delta^{\frac{1}{2}-2\epsilon} [B]_{-\epsilon} [v]_{3/2-\epsilon, \mcD_{s+s_1/4}},
$$
and using the estimate \eqref{eq:est_3_2} we get
\begin{equation*} \begin{split}
\big\Vert B_\delta\cdot\nabla v_\delta - (B\cdot\nabla v)_\delta \big\Vert_{\mcD_s+s_1/2} &\lesssim [B]_{-\epsilon} \Vert v\Vert_{\mcD_s} \Big(\frac{\delta}{s_1}\Big)^{\frac{1}{2}-2\epsilon} s_1^{-1+3\epsilon}   \\
&\lesssim [B]_{-\epsilon} \Vert v\Vert_{\mcD_s} k^{-\frac{1}{2}+2\epsilon} s_1^{-1+3\epsilon}.
\end{split} \end{equation*}
$\bullet$ We also use Lemma \ref{LemCommutatorEstimate} to deal with the $Z_2$ term and write
\begin{equation*} \begin{split}
(Z_2 v^2)_\delta &= (Z_2)_\delta v^2 + O\big(\delta^\epsilon[Z_2][v^2]_{1/2+2\epsilon,\mcD_{s+s_1/4}}\big)   \\
&= O\big(\delta^{-\frac{1}{2}-\epsilon}[Z_2]\Vert v\Vert_{\mcD_s}^2\big) + O\big(\delta^\epsilon [Z_2] \Vert v\Vert_{\mcD_s}[v]_{1/2+2\epsilon,\mcD_{s+s_1/4}}\big),
\end{split} \end{equation*}
and using the variation of \eqref{eq:est_1_2} we obtain the estimate on the $Z_2$ term in \eqref{EqBoundsStep2Z} since $\delta\leq s_1$. Similarly one has
\begin{equation*} \begin{split}
(Z_1 v)_\delta &= (Z_1)_\delta v + O\big(\delta^\epsilon[Z_1][v]_{2\epsilon,\mcD_{s+s_1/4}}\big)   \\
&= O\big(\delta^{-\frac{1}{2}-\epsilon}[Z_1] \Vert v\Vert_{\mcD_s}\big) + O\big(\delta^\epsilon [Z_1] \Vert v\Vert_{\mcD_s} [v]_{2\epsilon,\mcD_{s+s_1/4}}\big)   \\
&= O\big(\delta^{-\frac{1}{2}-\epsilon}[Z_1] \Vert v\Vert_{\mcD_s}\big) + O\big(\delta^\epsilon [Z_1] \Vert v\Vert_{\mcD_s} s_1^{-2\epsilon} ]\Vert v\Vert_{\mcD_s}\big)   \\
&= O\big(\delta^{-\frac{1}{2}-\epsilon}[Z_1] \Vert v\Vert_{\mcD_s}\big).
\end{split} \end{equation*}

\ssk

\noindent $\rhd$ \textbf{\textsf{Step 3. Choice of scales $\lambda_0$ and $\delta$.}} We choose
$$
s_1 \geq \lambda_0 = \Vert v\Vert_{\mcD_s}^{-1}, \qquad \delta = c_1 \Vert v\Vert_{\mcD_s}^{-1-\epsilon}
$$
for a positive constant $c_1\geq 1$ to be chosen below, and we take
\begin{equation} \label{EqChoiceK}
k = c_1^{-1} \Vert v\Vert_{\mcD_s}^\epsilon.
\end{equation}
One then has
\begin{equation} \label{EqLastEstimates1} \begin{split}
&\big\| [\mcL , (\cdot)_\delta](v) \big\|_{\mcD_{s+s_1/2}} \lesssim c_1^{-1} \Vert v\Vert_{\mcD_s}^{2+\epsilon},   \\
&\big\| A v_\delta^3 - (Av^3)_\delta \big\|_{\mcD_{s+s_1/2}} \lesssim c_1^{-\frac{1}{2}} \Vert A\Vert_{1/2,\mcD_s} \Vert v\Vert_{\mcD_s}^{3-\epsilon(\frac{1}{2}-\epsilon)},   \\
&\big\| B_\delta \cdot \nabla v_\delta - (B\cdot \nabla v )_\delta \big\|_{\mcD_{s+s_1/2}} \lesssim c_1^{\frac{1}{2}} [B]_{-\epsilon} \Vert v\Vert_{\mcD_s}^{2-\epsilon(\frac{1}{2}-2\epsilon)-3\epsilon},   \\
\end{split} \end{equation}
and
\begin{equation} \label{EqLastEstimates2} \begin{split}
&\| (Z_2 v^2)_\delta \|_{\mcD_{s+s_1/2}} \lesssim c_1^{-\frac{1}{2}-\epsilon} [Z_2]_{-1/2-\epsilon} \|v\|_{\mcD_s}^{2+(1+\epsilon)(\frac{1}{2}+\epsilon)},   \\
&\| (Z_1 v)_\delta \|_{\mcD_{s+s_1/2}} \lesssim c_1^{-\frac{1}{2}-\epsilon} [Z_1]_{-1/2-\epsilon} \Vert v\Vert_{\mcD_s}^{1+(1+\epsilon)(\frac{1}{2}+\epsilon)},   \\
&\| (Z_0)_\delta \|_{\mcD_{s+s_1/2}} \lesssim c_1^{-\frac{1}{2}-\epsilon} [Z_0]_{-1/2-\epsilon} \Vert v\Vert_{\mcD_s}^{(1+\epsilon)(\frac{1}{2}+\epsilon)}.
\end{split} \end{equation}
We then choose $c_1$ so that the estimate on $ [\mcL , (\cdot)_\delta](v)$ in \eqref{EqLastEstimates1} reads
$$
\big\| [\mcL , (\cdot)_\delta](v) \big\|_{\mcD_{s+s_1/2}} \leq 2^{-3} \, \Vert v\Vert_{\mcD_s}^{2+\epsilon}.
$$
Then one chooses $c$ in \eqref{EqAssumptionMWProof} so that the $A,B,Z_i$ terms in \eqref{EqLastEstimates1} and \eqref{EqLastEstimates2} are all smaller than $2^{-3} A_- \Vert v\Vert_{\mcD_s}^3$. The choice of exponents $m_\tau$ in \eqref{EqChoiceExponents}  was done precisely for that purpose. The estimates \eqref{EqEstimateStep2} and \eqref{EqVFromVdelta} together give
\begin{eqnarray}\label{eq:l_infty_main_est}
\| v\|_{\mcD_{s+s_1}} \leq \max \left\{ \frac{2A_-^{-1/2}}{s_1} \, , \, \frac{1}{2} \, \| v\|_{\mcD_s} \right\}.
\end{eqnarray}

\ssk

\noindent {\it (b) Proof of Theorem \ref{ThmLinftyComingDown} from \eqref{eq:l_infty_main_est}.} We now use an  argument due to Moinat \& Weber \cite{MoinatWeber} to derive the $L^\infty$ bound \eqref{EqMWBound} from the estimate \eqref{eq:l_infty_main_est}. We proceed differently depending on whether $\min_{\mcD_s} A \leq 1$ or $\min_{\mcD_s} A > 1$. For $\min_{\mcD_s} A \leq 1$ we have for
$$
\widetilde v \defeq A_-^{1/2} \, v.
$$ 
the estimate
\begin{equation}\label{eq:rewrite_est}
\|\widetilde v\|_{\mcD_{s+s'}} \leq \max \left\{ \frac{2}{s'}, \frac{1}{2} \, \|\widetilde v\|_{\mcD_s} \right\},
\end{equation}
We set $s_1\defeq 4\,\Vert\widetilde{v}\Vert_{\mcD_s}^{-1}$ and define the times $s=s_0<s+s_1<\ldots<s+s_N=\frac{1}{2}$ from the relation 
\begin{equation}\label{eq:sequence}
s_{n+1}-s_{n} = 4\,\|\widetilde v\|_{\mcD_{s+s_n}}^{-1}.
\end{equation}
The sequence terminates once $s+s_{n+1}\geq 1/2$, in which case we set $s_{n+1}=s_{N}= 1/2$, or once the assumption \eqref{EqAssumptionMWProof} fails for $\mcD_{s+s_{n+1}}$. Since $4\|\widetilde v\|_{\mcD_{s+s_n}}^{-1}$ is increasing in $n$ the sequence terminates after finitely many steps and we have
\begin{eqnarray}\label{eq:reduce_1/2}
\|\widetilde v\|_{\mcD_{s+s_{n}}} \leq \frac{1}{2}\|\widetilde v\|_{\mcD_{s+s_{n-1}}}.
\end{eqnarray}
It follows from \eqref{eq:sequence} and \eqref{eq:reduce_1/2} that for $1 \leq n \leq N-1$
\begin{eqnarray}\label{eq:est_s_n}
s+s_n \lesssim \| \widetilde v\|_{\mcD_{s+s_n}}^{-1}.
\end{eqnarray} 
Therefore  the bound in  Theorem \ref{ThmLinftyComingDown} holds at the time $s+s_n$, for $n\leq N-1$. For the last domain $\mcD_{s_N}$, if the assumption \eqref{EqAssumptionMWProof} fails for $\mcD_{s_N}$, then we get the bound immediately. Otherwise we have either $s+s_{N-1}\geq 1/4$ or $s_N-s_{N-1}\geq 1/4$. In the first case, we use the  estimate \eqref{eq:est_s_n} for $s+s_{N-1} $
\begin{eqnarray}\label{eq_s_N_1}
s+s_{N} = \frac{1}{2}\leq 2(s+ s_{N-1})\lesssim \|\widetilde v\|_{\mcD_{s+s_{N-1}}}^{-1}.
\end{eqnarray}
In the second case, by \eqref{eq:sequence}, we have
\begin{equation}\label{eq_s_N_2}
s+s_N = \frac{1}{2}\leq 2 (s_{N}-s_{N-1}) = 4 \, \|\widetilde v\|_{\mcD_{s+s_{N-1}}}^{-1}.
\end{equation}
Finally, for any time $r\in (s+s_n, s+s_{n+1})$ with $0\leq n\leq N-2$, using \eqref{eq:sequence} and \eqref{eq:est_s_n} we infer that
$$
r\leq s+s_{n+1} = s+s_n + (s_{n+1}-s_n) \lesssim \|\widetilde v\|_{\mcD_{s+s_n}}^{-1} \lesssim \|\widetilde v\|_{\mcD_r}^{-1},
$$
and for $t\in (s+s_{N-1}, s+s_N)$, \eqref{eq_s_N_1} and \eqref{eq_s_N_2} imply that
$$
r\leq s+s_N \lesssim \| \widetilde v\|_{\mcD_{s+s_{N-1}}}^{-1} \leq \| \widetilde v\|_{\mcD_r}^{-1}.
$$
This gives the desired estimate
$$
\Vert v\Vert_{\mcD_s} \leq \frac{2A_-^{-\frac{1}{2}}}{s} .
$$
In the case where $\min_{\mcD_s} A > 1$ we infer from \eqref{eq:l_infty_main_est} that
$$
\| v\|_{\mcD_{s+s_1}} \leq \max \left\{ \frac{2}{s_1}, \frac{1}{2} \| v\|_{\mcD_s} \right\}.
$$
We get the estimate $\| v\|_{\mcD_{s}}\leq 2/s$ by repeating the preceding argument.

\bigskip

\section{Controlling stronger norms of $v$}
\label{SubsectionControlNorm32}

The next statement requires that $0<\epsilon\leq 1/4$. As a shorthand notation we write below $\Vert\tau\Vert=\Vert\tau\Vert_{C([0,T],C^{\vert\tau\vert}(M))}$ for $\tau\in\{A,B,Z_2,Z_1\}$.

\medskip

\begin{thm} \label{ThmHighControl}
There are two functions $C_\emptyset$ and $C'_\emptyset$ of $\Vert A\Vert, \Vert B\Vert, \Vert Z_2\Vert, \Vert Z_1\Vert$ that do \emph{not} depend on $Z_0$ such that setting
$$
\lambda_0 = C_\emptyset \wedge \big(\max(\Vert Z_1\Vert,\Vert Z_2\Vert) \Vert v\Vert_{\mcD_s}\big)^{-\frac{2}{3-4\epsilon}}   \vspace{-0.1cm}
$$
one has
$$
[v]_{3/2-\epsilon, \mcD_{s+\lambda_0}} \lesssim C'_\emptyset \big(1\vee \Vert v\Vert_{\mcD_s}\big)^3 + \vert Z_0\vert.
$$
\end{thm}

\ssk

The precise formulas for $C_\emptyset$ and $C'_\emptyset$ can be extracted from the proof of Theorem \ref{ThmHighControl} and have no importance in this work.

\ssk

\begin{Dem}
The proof of this statement is essentially a variation on the content of Step 1 in the proof of Theorem \ref{ThmLinftyComingDown}. We repeat it here for the reader's convenience. We start with the equation

\begin{equation} \label{LVDelta}
(\mcL v)_\delta = - (Av^3)_\delta + (B\cdot \nabla v)_\delta + (Z_2 v^2)_\delta + (Z_1v)_\delta + (Z_0)_\delta. 
\end{equation}
Note that the conclusion of Proposition \ref{LemLocalSchauder} makes it possible to use the $(3/2-\epsilon)$-H\"older seminorm of $v$ in our estimates of the right-hand side of \eqref{LVDelta} provided it comes with a small factor that can eventually be absorbed in the left-hand side of \eqref{EqSchauderMain}. We work with $0<\delta\leq \lambda/4$. For the $Z_0$ and $a$ terms in \eqref{LVDelta} we have
\begin{equation*} \begin{split}
\| (Z_0)_\delta\|_{\mcD_{s+\lambda-\delta}} &\leq \delta^{-\frac{1}{2} - \epsilon}\Vert Z_0\Vert   \\
\| (Av^3)_\delta \|_{\mcD_{s+\lambda-\delta}} &\leq \|A\| \|v\|_{\mcD_s}^{3}. 
\end{split} \end{equation*}
For the $Z_2$ and $Z_1$ terms in \eqref{LVDelta} one gets

\begin{equation*} \begin{split}
\|(Z_2v^2)_{\delta}\|_{\mcD_{s+\lambda-\delta}} &\leq \delta^{-\frac{1}{2}-\epsilon} \Vert Z_2\Vert \,  \|v\|_{\mcD_s}^{2} + 2  \delta^{-\epsilon} \Vert Z_2\Vert \, [v]_{(1/2+2\epsilon\,;\,\delta), \mcD_{s +\lambda/2}}   \|v\|_{\mcD_s},   \\
\|(Z_1v)_{\delta}\|_{\mcD_{s+\lambda- \delta}} &\leq \delta^{-\frac{1}{2}-\epsilon |} \Vert Z_1\Vert \, \|v\|_{\mcD_s} + \delta^{-\epsilon} \Vert Z_1\Vert \, [v]_{(1/2+2\epsilon\,;\,\delta), \mcD_{s+\lambda/2}}.
\end{split} \end{equation*}
It turns out to be useful to introduce the commutator $B_\delta\cdot\nabla v - (B\cdot\nabla v)_\delta$ to estimate $(B\cdot\nabla v)_\delta$ itself as we get from Lemma \ref{LemCommutatorEstimate} the bound
\begin{equation*} \begin{split}
\Vert (B\cdot\nabla v)_\delta\Vert_{\mcD_{s+\lambda-\delta}} &\leq \big\Vert B_\delta\cdot\nabla v - (B\cdot\nabla v)_\delta \big\Vert_{\mcD_{s+\lambda-\delta}} + \Vert B_\delta\cdot\nabla v\Vert_{\mcD_{s+\lambda-\delta}}   \\
&\lesssim \delta^{\frac{1}{2}-2\epsilon} \Vert B\Vert \, [\nabla v]_{1/2-\epsilon, \mcD_{s+\lambda/2}} + \delta^{-\epsilon} \Vert B\Vert\, \Vert \nabla v\Vert_{\mcD_{s+\lambda/2}}.
\end{split} \end{equation*}
For all $\lambda$ sufficiently small one has
\begin{equation} \label{EqGradientBound} \begin{split}
\Vert \nabla v\Vert_{\mcD_{s+\lambda/2}} &\lesssim \lambda^{\frac{1}{2}-\epsilon} [v]_{3/2-\epsilon, \mcD_{s+\lambda/2}} + \lambda^{-1} \Vert v\Vert_{(0;\lambda), \mcD_{s+\lambda/2}},   \\
[\nabla v]_{1/2-\epsilon, \mcD_{s+\lambda/2}} &\lesssim [v]_{3/2-\epsilon, \mcD_{s+\lambda/2}} + \lambda^{-\frac{3}{2}+\epsilon} \Vert v\Vert_{(0;\lambda), \mcD_{s+\lambda/2}}.
\end{split} \end{equation}
Combining these estimates all together and $\|v\|_{(0;\lambda), \mcD_{s+\lambda/2}} \leq 2\|v\|_{\mcD_s}$ yields the bound
$$
\delta^{2-\frac{3}{2}+\epsilon} \|(\mcL v)_\delta \|_{\mcD _{s+\lambda-\delta}} \lesssim C_\lambda   \vspace{-0.1cm}
$$
with 
\begin{equation*} \begin{split}
C_\lambda \defeq \lambda^{\frac{1}{2}+\epsilon} &\|A\| \|v\|_{\mcD_s}^{3} + \Vert B\Vert \, [v]_{3/2-\epsilon, \mcD_{s+\lambda/2}} + \lambda^{-\frac{1}{2}} \Vert B\Vert \, \|v\|_{\mcD_s}   \\
&+ \Vert Z_2\Vert_{-1/2-\epsilon} \|v\|_{\mcD_s}^{2} +  \lambda^{\frac{1}{2}}\Vert Z_2\Vert _{-1/2-\epsilon}[v]_{(1/2+2\epsilon\,;\,\delta), \mcD_{s+\lambda/2}}   \|v\|_{\mcD_s}  \\
&+ \Vert Z_1\Vert \, \|v\|_{\mcD_s} + \lambda^{\frac{1}{2}}\Vert Z_1\Vert \, [v]_{(1/2+2\epsilon\,;\,\delta), \mcD_{s+\lambda/2}}   \|v\|_{\mcD_s} + \Vert Z_0\Vert.
\end{split} \end{equation*}
The Schauder estimate from Proposition \ref{LemLocalSchauder} gives us
\begin{equation} \label{EqAlmostBound32Norm}
\underset{\lambda\leq\frac{\lambda_0}{2}}{\sup}\; \lambda^{\frac{3}{2}-\epsilon} [u]_{3/2-\epsilon,\mcD_{s+\lambda}} \leq \underset{\lambda\leq\lambda_0}{\sup}\; \lambda^{\frac{3}{2}-\epsilon} [u]_{3/2-\epsilon,\mcD_{s+\lambda}} \lesssim \underset{\lambda\leq\lambda_0/2}{\sup}\; \lambda^{\frac{3}{2}-\epsilon} C_\lambda + \Vert v\Vert_{\mcD_s}
\end{equation}
Even though the $C_\lambda$ depend on the $(3/2-\epsilon)$-seminorm of $v$ on $\mcD_{s+\lambda}$ it comes with a small multiplicative factor if one chooses
\begin{equation} \label{EqDefnLambada0}
\lambda_0 = \Big(4 \max\Big\{ \Vert B\Vert^{\frac{1}{1-\epsilon}} \, ; \, \Vert Z_1\Vert^{\frac{1}{3/2-2\epsilon}} \, ;\, \big(\Vert v\Vert_{\mcD_{s}} \Vert Z_2\Vert \big)^{\frac{1}{3/2-2\epsilon}} \Big\}\Big)^{-1}.
\end{equation}
The term of $C_\lambda$ that involves the $(3/2-\epsilon)$-seminorm of $v$ can then be absorbed in the left-hand side of \eqref{EqAlmostBound32Norm}. Still the constant $C_\lambda$ depends on some $(1/2+2\epsilon)$-seminorm of $v$. To eventually have a bound on $C_\lambda$ that only involves $\|v\|_{ \mcD_s}$ we use the elementary estimate
\begin{equation*}
[v]_{(1/2+2\epsilon\,;\,\delta), \mcD_{s+\lambda/2}} \leq [v]_{1/2+2\epsilon,\mcD_{s+\lambda/2}} \leq \lambda^{1-3\epsilon} [v]_{3/2-\epsilon,\mcD_{s+\lambda/2}} + \lambda^{\frac{1}{2}-2\epsilon} \Vert \nabla v\Vert_{\mcD_{s+\lambda/2}}
\end{equation*}
and the gradient bound \eqref{EqGradientBound} on $v$ in uniform norm to see that 
$$
[v]_{(1/2+2\epsilon;\delta), \mcD_{s+\lambda/2}} \lesssim \lambda^{1-3\epsilon} [v]_{3/2-\epsilon,\mcD_{s+\lambda/2}} + \lambda^{-\frac{1}{2}-2\epsilon} \Vert v\Vert_{\mcD_{s+\lambda/2}}.
$$
One therefore has for $0<\epsilon\leq \epsilon_0$
$$
\lambda^{2-\epsilon_0+2\epsilon}\Vert Z_i\Vert \, [v]_{(1/2-\epsilon\,;\,\delta), \mcD_{s+\lambda/2}} \lesssim \Vert Z_i\Vert  \big( \lambda^{3-\epsilon_0-2\epsilon} [v]_{ 3/2-\epsilon,\mcD_{s+\lambda/2}}+ \lambda^{\frac{3}{2}-\epsilon_0} \Vert v\Vert_{\mcD_{s}}\big)
$$ 
for $1\leq i\leq 2$, so our choice of $\lambda_0$ ensures that this term can also be absorbed in the left-hand side of \eqref{EqAlmostBound32Norm}. We get in the end the bound
$$
\sup_{\lambda\leq \lambda_0/2} \lambda^{\frac{3}{2}- \epsilon} [v]_{\frac{3}{2}-\epsilon, \mcD_{s+\lambda}} \lesssim \sup_{\lambda\leq \lambda_0} C'_\lambda + \Vert v\Vert_{\mcD_{s}}
$$
where 
\begin{equation*} \begin{split}
C'_\lambda =  \lambda^2 \|A\| \, \|v\|_{\mcD_s}^3 &+ \lambda^{1-\epsilon} \Vert B\Vert \, \|v\|_{ \mcD_s} + \lambda^{\frac{3}{2} -\epsilon} \Vert Z_2\Vert \, \|v\|_{ \mcD_s}^{2} +  \lambda^{\frac{3}{2}-\epsilon_0} \Vert Z_1\Vert \, \|v\|_{ \mcD_s}   \\
&+ \lambda^{\frac{3}{2}}  \Vert Z_2\Vert \, \|v\|^2_{ \mcD_s}  + \lambda^{\frac{3}{2}} \Vert Z_1\Vert \, \|v\|_{\mcD_s} + \lambda^{\frac{3}{2}-\epsilon}\Vert Z_0\Vert.
\end{split} \end{equation*}
This is an increasing function of $\lambda$ and writing 
$$
[v]_{\frac{3}{2}-\epsilon, \mcD_{s+\lambda_0}} \lesssim \lambda_0^{-\frac{3}{2}+ \epsilon} \big(C'_{\lambda_0} + \Vert v\Vert_{\mcD_{s}}\big)
$$
gives the conclusion in view of the definition \eqref{EqDefnLambada0} of $\lambda_0$. 
\end{Dem}

\medskip

\section{Variation on a theme}
\label{SectionVariationTheme}

We used in Section \ref{SectionUniqueness} a perturbed version of the dynamics \eqref{EqJPFormulationSynthetic} that contains an additional drift with a particular form. This section is dedicated to proving for the solution of this perturbed dynamics some explicit control on its $L^\infty$ and stronger norm similar to Theorem \ref{ThmLinftyComingDown} and Theorem \ref{ThmHighControl}. We use the convention that ${\bf 1}_{1<\cdot<\tau}=0$ if $\tau=1$. Below, the space $\llparenthesis\alpha_0,1+2\epsilon\rrparenthesis$ was introduced in Section 2.1.2 of \cite{BDFT} for a particular value of $\alpha_0$ that does not matter here. The norm on this space quantifies the explosion of a function of time $t>0$ as $t$ goes to $0$. Its precise definition here does not really matter and we only note that for a fixed $0<T<\tau$ the restriction to $[T,\tau)$ of an element of the space $C([0,\tau),C^{-1/2-\epsilon}(M)) \cap \llparenthesis \alpha_0, 1+2\epsilon\rrparenthesis$ is an element of $C([T,\tau],C^{1+2\epsilon}(M))$. Let $v$ stand below for the solution of \eqref{EqJPFormulationSynthetic} started from some initial condition $\phi_1'$.

\ssk

\begin{thm} \label{ThmVariationOnTheme}
Pick a constant $\ell\in\bbR$ and an initial condition $\phi_2' \in C^{-1/2-\epsilon}(M)$. The equation
\begin{equation} \label{EqModifiedJPEquation} \begin{split}
\partial_tv_\ell &= (\Delta-1)v_\ell - Av_\ell^3 + B\nabla v_\ell + Z_2v_\ell^2 + Z_1v_\ell + Z_0   \\
&\qquad+ \ell \, {\bf 1}_{1< t <\tau}\,\frac{v(t)-v_\ell(t)}{\Vert v(t)-v_\ell(t)\Vert_{L^2}}\,\exp\big(3\begin{tikzpicture}[scale=0.3,baseline=0cm] \node at (0,0) [dot] (0) {}; \node at (0,0.4) [dot] (1) {}; \node at (-0.3,0.8)  [noise] (noise1) {}; \node at (0.3,0.8)  [noise] (noise2) {}; \draw[K] (0) to (1); \draw[K] (1) to (noise1); \draw[K] (1) to (noise2); \end{tikzpicture}(t)\big)   \\
&\eqdef F(v_\ell) + \ell \, {\bf 1}_{1< t <\tau}\,\frac{v(t)-v_\ell(t)}{\Vert v(t)-v_\ell(t)\Vert_{L^2}} \, \exp\big(3\begin{tikzpicture}[scale=0.3,baseline=0cm] \node at (0,0) [dot] (0) {}; \node at (0,0.4) [dot] (1) {}; \node at (-0.3,0.8)  [noise] (noise1) {}; \node at (0.3,0.8)  [noise] (noise2) {}; \draw[K] (0) to (1); \draw[K] (1) to (noise1); \draw[K] (1) to (noise2); \end{tikzpicture}(t)\big) \eqdef F_\ell(t,v,v_\ell), \qquad (0\leq t < \tau),
\end{split} \end{equation}
where 
$$
\tau = \tau(\ell,\phi_1,\phi_2) \defeq \inf\big\{s\geq 1 \, ; \, v_\ell(s) = v(s)\big\} \wedge 2,
$$
has a unique solution in $C\big([0,\tau),C^{-1/2-\epsilon}(M)\big) \cap \llparenthesis \alpha_0, 1+2\epsilon\rrparenthesis$. Furthermore it satisfies the estimates
\begin{equation} \label{EqControlVEllInfty}
\Vert v_\ell(s)\Vert_{L^\infty} \leq c_1(\widehat{\xi}\,) \, (1 + {\color{black} \ell^{\frac{1}{3}}})
\end{equation}
and 
\begin{equation} \label{EqControlVEll1Plus}
{\color{black} \Vert v_\ell(s)\Vert_{C^{1+2\epsilon}}} \leq c_2(\widehat{\xi}\,) \, (1 + {\color{black} \ell})
\end{equation}
for all $1\leq s<\tau$, for some explicit functions $c_1(\widehat{\xi}\,), c_2(\widehat{\xi}\,)$ of $\widehat{\xi}$ whose precise values play no role in what follows.
\end{thm}

\ssk

\begin{Dem}
{\it Local in time well-posedness beyond time $1$.} There is no loss of generality in assuming that $v(1,\phi_1')\neq v_\ell(1,\phi_2')$. Denote by $C_{v_\ell(1,\phi_2')}\big([1,T], C^{1+2\epsilon}(M)\big)$ the set of continuous functions from the interval $[1,T]$ into $C^{1+2\epsilon}(M)$ with value $v_\ell(1,\phi_2')$ at time $1$. Pick a positive constant 
$$
m < \big\Vert v(1,\phi_1') - v_\ell(1,\phi_2')\big\Vert_{L^2} \wedge 1
$$
(think of it as being small) and set 
$$
\mcV_{v_\ell(1,\phi_2')}(m,T) \defeq \Big\{w\in C_{v_\ell(1,\phi_2')}\big([1,T], C^{1+2\epsilon}(M)\big) \, ; \, \min_{1\leq t\leq T}\, \Vert v(t) - w(t)\Vert_{L^2} > m \Big\}.
$$
This is an open subset of $C_{v_\ell(1,\phi_2')}([1,T], C^{1+2\epsilon}(M))$ with closure $\overline\mcV_{v_\ell(1,\phi_2')}(m,T)$ included in
$$
\Big\{w\in C_{v_\ell(1,\phi_2')}\big([1,T], C^{1+2\epsilon}(M)\big)\,;\, \min_{1\leq t\leq T}\, \Vert v(t) - w(t)\Vert_{L^2} \geq m\Big\}.
$$ 

\ssk

\begin{lem} \label{LemLocalInTimeWellPosedness}
There exists a positive time $T(\ell,m)$ such that the map $\mathscr{F}$ defined as
$$
\mathscr{F}(w)(t) \defeq e^{(t-1)(\Delta-1)} \big(v_\ell(1,\phi'_2)\big) + \int_0^{t-1} e^{(t-1-s)(\Delta-1)}\big(F_\ell(1+s,w)\big)ds
$$
is a contraction of $\overline\mcV_{v_\ell(1,\phi_2')}\big(m,T(\ell,m)\big)$ into itself. One can choose $T(\ell,m)$ as a decreasing function of $m$.
\end{lem}

\ssk

\textbf{\textsf{Proof --}} Indeed, for $w\in \mcV_{v_\ell(1,\phi_2')}(C,T)$ one has $F_\ell(w)\in C\big([1,T],C^{-1/2-\epsilon}(M)\big)$, and for any $w_1,w_2$ in $\mcV_{v_\ell(1,\phi_2')}(C,T)$ one has 
\begin{equation*} \begin{split}
F_\ell(w_1) - F_\ell(w_2) &= (\Delta-1)(w_1-w_2) - A \big(w_1^3-w_2^3\big) + B\nabla(w_1-w_2)   \\
&\quad+ Z_2 \big(w_1^2-w_2^2\big) + Z_1(w_1-w_2) + \ell\Big(\frac{v-w_1}{\Vert v-w_1\Vert_{L^2}} - \frac{v-w_2}{\Vert v-w_2\Vert_{L^2}}\Big)
\end{split} \end{equation*}
and
\begin{equation*} \begin{split}
\big\Vert F_\ell(w_1) - F_\ell(w_2) \big\Vert_{C_TC^{-1/2-\epsilon}} &\lesssim 1 + 2\big(\Vert w_1\Vert_{L^\infty}^2 + \Vert w_2\Vert_{L^\infty}^2\big) \Vert w_1-w_2\Vert_{L^\infty} + \Vert w_1-w_2\Vert_{C^{1+\epsilon}}   \\
&\quad+ \big(\Vert w_1\Vert_{L^\infty} + \Vert w_2\Vert_{L^\infty} + 1\big) \Vert w_1-w_2\Vert_{C^{1/2+\epsilon}}   \\
&\quad+ \ell \, m^{-2}\big(\Vert v\Vert_{L^\infty} + \Vert w_1\Vert_{L^\infty} + \Vert w_2\Vert_{L^\infty}\big) \Vert w_1-w_2\Vert_{L^\infty}.
\end{split} \end{equation*}
(We estimated the $\ell$-term in $L^\infty$ by bounding the operator norm of the map $\Vert f\Vert_{L^2}f$ from $L^\infty$ into itself. We left the constants on the right hand side to make the verification easier.) The small-time estimate
$$
\big\Vert \big(\partial_t-(\Delta-1)\big)^{-1}(f)\big\Vert_{C_TC^{1+2\epsilon}} \lesssim T^{\frac{1}{2}-3\epsilon} \Vert f\Vert_{C_TC^{-1/2-\epsilon}},
$$
then entails the statement of the lemma.   \hfill $\RHD$

\ssk

Two solutions corresponding to $m_2 < m_1$ satisfy $T(\ell,m_2)\geq T(\ell,m_1)$ and coincide on the interval $[0,T(\ell,m_1)]$ from uniqueness. We prove the non-explosion of the solution to Equation \eqref{EqModifiedJPEquation} before the coupling time with the path $(u(t))_{t\geq 0}$ by showing that any solution to \eqref{EqModifiedJPEquation} on a time interval $(0,T]$ satisfies the size estimates \eqref{EqControlVEllInfty} and \eqref{EqControlVEll1Plus}. By repeated use of the local in time result of Lemma \ref{LemLocalInTimeWellPosedness} these estimates allow to define the maximal existence time $\tau(\ell,m)\geq T(\ell,m)$ before $\Vert v_\ell(t)-v(t)\Vert_{L^2}=m$, with $\tau(\ell,m_2)\geq \tau(\ell,m_1)$ if $m_2 < m_1$. The coupling time $\tau$ is defined as
$$
\tau = \tau(\ell) = \underset{m\downarrow 0}{\lim}\, \tau(\ell,m).
$$
We assume in the next paragraph that $v_\ell$ is a solution to \eqref{EqModifiedJPEquation} on a time interval $(0,T]$, so $\Vert v(t)-v_\ell(t)\Vert_{L^2}\geq m$ for some $m>0$ below.

\ssk

{\it Quantitative estimates \eqref{EqControlVEllInfty} and \eqref{EqControlVEll1Plus}.} The proof of these statements is similar to the proofs of Theorem \ref{ThmLinftyComingDown} and Theorem \ref{ThmHighControl} but one has to deal with the fact that the drift
$$
G_\ell(v,v_\ell)(t) \defeq \ell \,\frac{v(t)-v_\ell(t)}{\Vert v(t)-v_\ell(t)\Vert_{L^2}}\,\exp\big(3\begin{tikzpicture}[scale=0.3,baseline=0cm] \node at (0,0) [dot] (0) {}; \node at (0,0.4) [dot] (1) {}; \node at (-0.3,0.8)  [noise] (noise1) {}; \node at (0.3,0.8)  [noise] (noise2) {}; \draw[K] (0) to (1); \draw[K] (1) to (noise1); \draw[K] (1) to (noise2); \end{tikzpicture}(t)\big)
$$
does not have good controls in $C^{-1/2-\epsilon}(M)$ for all times. Rather it has good controls when seen as an element of $L^2(M)$. We only point out the modifications of the proof of Theorem \ref{ThmLinftyComingDown} needed to provide a proof of \eqref{EqControlVEllInfty} and \eqref{EqControlVEll1Plus}. We go linearly along the proof of the former.

\ssk

$\bullet$ In Step 1, looking at the drift as a continuous function of time with values in a fixed ball of $L^2(M)$ of radius proportional to $\ell$, in addition to the terms that were already analyzed in the Step 1 of the proof of Theorem \ref{ThmLinftyComingDown} one gets here an additional term
$$
\big\Vert \big(G_\ell(v,v_\ell)\big)_\delta\big\Vert_{\mcD_{s+\lambda-\delta}} \lesssim_{\widehat\xi} \delta^{-\frac{3}{4}} \, \ell.
$$
The worst exponent in Step 1 in Section {\sf \ref{SectionLinftyComingDown}} was $-1/2-\epsilon$. This difference between the worst exponents in the two situations explains why the reasoning of Step 1 cannot give control of the $3/2-\epsilon$ seminorm of $v_\ell$ as the quantity $\delta^{1/2+\epsilon} \|(\mcL v_\ell)_\delta \|_{\mcD _{s+\lambda-\delta}}$ is not bounded anymore. It gives however a control of the $(1+2\epsilon)$-seminorm of $v_\ell$ in terms of $\Vert v_\ell\Vert_{\mcD_s} + \ell$, following verbatim what was done in Section {\sf \ref{SectionProofLInftyComingDown}}.

$\bullet$ The contribution of the drift $G_\ell(v,v_\ell)$ in Step 2 is the same as in Step 1 and one needs in addition to replace here the use of the $(3/2-\epsilon)$-seminorm of $v$ made in Section {\sf \ref{SectionProofLInftyComingDown}} in the estimate of the $B$ term by the use of the $(1+2\epsilon)$-seminorm of $v_\ell$. 

$\bullet$ Step 3 works similarly as in Section {\sf \ref{SectionProofLInftyComingDown}} and provides the estimate
$$
\| v_\ell\|_{\mcD_{s+s_1}} \leq \max \left\{ \frac{2A_-^{-\frac{1}{2}}}{s_1} \, , \, \frac{1}{2} \, \| v_\ell\|_{\mcD_s} + \ell^{\frac{1}{3}}\right\},
$$
so we have
$$
\| v_\ell\|_{\mcD_{s+s_1}} \leq \max \left\{ \frac{2A_-^{-\frac{1}{2}}}{s_1} \, , \, \frac{3}{4} \, \| v_\ell\|_{\mcD_s}\right\}
$$
as long as $\Vert v_\ell\Vert_{\mcD_s} \geq 4\,\ell^{\frac{1}{3}}$. Proceeding as in part {\it (b)} of the proof of Theorem \ref{ThmLinftyComingDown} with the contraction coefficient $1/2$ replaced by $3/4$, one gets
$$
\Vert v\Vert_{\mcD_s} \lesssim \max\bigg\{\frac{1}{s} , \big((c_A [\tau]_{\vert\tau\vert})^{m_\tau}\big)_{\tau\in\mcT} \, , \, 4\,\ell^{\frac{1}{3}}\bigg\}.
$$
The estimate \eqref{EqControlVEllInfty} follows from the preceding inequality. We repeat the proof of Theorem \ref{ThmHighControl} to obtain \eqref{EqControlVEll1Plus}, adding the contribution of the drift $G_\ell(v,v_\ell)$ and replacing the $(3/2-\epsilon)$-seminorm by the $(1+2\epsilon)$-seminorm.
\end{Dem}

\medskip

\appendix
\section{An elementary continuity result}
\label{SectionAppendix}

The following statement is a simple variation on the classical proof of continuity of the paraproduct and resonant operators; we learned it from V.N. Dang although it is likely to be known already.

\ssk

\begin{lem} \label{LemSmallFactorTrick}
Let $(p_1,q_1), (p_2,q_2), (p,q)$ in $[1,+\infty]$ be such that 
$$
\frac{1}{p_1} + \frac{1}{p_2} = \frac{1}{p},
$$
$$
\frac{1}{q_1} + \frac{1}{q_2} = \frac{1}{q}.
$$
For any $\gamma>0$ there is a constant $C_\gamma$ such that for all integers $N$ and real numbers $a_1\geq 0$ one has
	\begin{equation} \label{EqTrickParaproduct}
	\Vert f\prec g\Vert_{B^{a_2-\gamma}_{pq}} \leq C_\gamma \, \Vert f\Vert_{B^{a_1}_{p_1q_1}} \Big( 2^{-N\gamma} \Vert g\Vert_{B^{a_2}_{p_2q_2}} + N \Vert g\Vert_{L^{p_2}}\Big).
	\end{equation}
\end{lem}

\ssk

We give here the details for the reader's convenience, when things are set in a Euclidean space. An elementary adaptation of the pattern of proof is needed to make it work in the setting of a $3$-dimensional Riemannian manifold with the paraproduct and resonant operators defined in Appendix {\sf A} of \cite{BDFT} or in \cite{BDFTCompanion}. We use the usual $\Delta_k$ notation for the Littlewood-Paley projectors.

\ssk

\begin{Dem}
Write
\begin{eqnarray}
f \prec g = \sum_{\ell< N, k\leq \ell-2} (\Delta_k f) (\Delta_\ell f) + \sum_{\ell\geq N,  k\leq \ell-2} (\Delta_k f) (\Delta_\ell g)
\end{eqnarray}
for any integer $N$. On the one hand, one has for all $m\in\bbN$

\begin{equation*} \begin{split}
\bigg\Vert \Delta_m \bigg(\sum_{\ell < N, k\leq \ell-2} (\Delta_k f)(\Delta_\ell g)\bigg) \bigg\Vert_{L^p} &\lesssim {\bf 1}_{m-1\leq N} \sum_{\vert\ell-m\vert\leq 1, k\leq \ell-2} \big\Vert (\Delta_k f)(\Delta_\ell g)\big\Vert_{L^p}   \\
\end{split} \end{equation*}
\begin{equation*} \begin{split}
&\lesssim  {\bf 1}_{m-1\leq N} \sum_{\vert\ell-m\vert\leq 1, k\leq \ell-2} 2^{-\ell a_1} \, 2^{ka_1}\Vert\Delta_k f\Vert_{L^{p_1}} \Vert \Delta_\ell g\Vert_{L^{p_2}}   \\
&\lesssim  {\bf 1}_{m-1\leq N} \, m\, \Vert f\Vert_{B^{a_1}_{p_1,q_1}}\Vert g\Vert_{L^{p_2}} \lesssim N \Vert f\Vert_{B^{a_1}_{p_1,q_1}}\Vert g\Vert_{L^{p_2}},
\end{split} \end{equation*}
using in the penultimate inequality H\"older inequality and the fact that $a_1\geq 0$, and Young convolution inequality to see that $\Vert \Delta_\ell g\Vert_{L^{p_2}} \leq \Vert g\Vert_{L^{p_2}}$. On the other hand, one has for $m\geq N-1$

\begin{equation*} \begin{split}
\bigg\Vert \Delta_m \bigg(\sum_{\ell\geq N, k\leq \ell-2} (\Delta_k f)(\Delta_\ell g)\bigg) \bigg\Vert_{L^p} &\lesssim \sum_{\vert\ell-m\vert\leq 1, k\leq \ell-2} \big\Vert (\Delta_k f)(\Delta_\ell g)\big\Vert_{L^p}   \\
&\lesssim \sum_{\vert\ell-m\vert\leq 1, k\leq \ell-2} 2^{-k a_1} \, 2^{ka_1}\Vert\Delta_k f\Vert_{L^{p_1}} \Vert \Delta_\ell g\Vert_{L^{p_2}}   \\
\end{split} \end{equation*}
\begin{equation*} \begin{split}
&\lesssim m^{\frac{q_1-1}{q_1}}\, \Vert f\Vert_{B^{a_1}_{p_1,q_1}} \sum_{\vert\ell-m\vert\leq 1} \Vert \Delta_\ell g\Vert_{L^{p_2}}
\end{split} \end{equation*}
from H\"older inequality in the sum over $k$. Estimate \eqref{EqTrickParaproduct} follows as a consequence. 
\end{Dem}

\end{document}